\documentclass[a4paper,12pt]{amsart}
\pdfoutput=1
\usepackage{preprint}   
\usepackage{multicol}
\usepackage{amsthm}
\usepackage{amsfonts}
\usepackage{graphicx}
\usepackage{verbatim}

\usepackage{subfig}

\graphicspath{{PDF/}{EPS/}}

\swapnumbers

\usepackage{amsmath}

\usepackage{fancybox}

\nonstopmode

\title[Computation of exterior moduli]
{Computation of exterior moduli of quadrilaterals}

\author[H.Hakula]{\noindent Harri Hakula}
\email{harri.hakula@tkk.fi}
\address{Aalto University, Institute of Mathematics,
P.O. Box 11100, FI-00076 Aalto,
FINLAND}

\author[A.Rasila]{\noindent Antti Rasila}
\email{antti.rasila@iki.fi}
\address{Aalto University, Institute of Mathematics,
P.O. Box 11100, FI-00076 Aalto,
FINLAND}

\author[M.Vuorinen]{\noindent Matti Vuorinen}
\email{vuorinen@utu.fi}
\address{Department of Mathematics,
FI-20014 University of Turku,
FINLAND}

\newcounter{minutes}
\divide\time by 60
\newcounter{hours}
\multiply\time by 60
\addtocounter{minutes}{-\time}

\keywords{conformal capacity, conformal modulus, quadrilateral modulus, $hp$-FEM, numerical conformal mapping}
\subjclass{65E05, 31A15, 30C85}

\begin{document}

\def\thefootnote{}
\footnotetext{ \texttt{\tiny File:~\jobname .tex,
         printed: \number\year-\number\month-\number\day,
         \thehours.\ifnum\theminutes<10{0}\fi\theminutes}
} \makeatletter\def\thefootnote{\@arabic\c@footnote}\makeatother

\begin{abstract}
We study the problem of computing the exterior modulus of a bounded
quadrilateral. We reduce this problem to the numerical solution of the
Dirichlet-Neumann problem for the Laplace equation.
Several experimental results, with error estimates, are
reported. Our main method makes use of an $hp$-FEM algorithm, which enables computations in the case of complicated geometry. For simple geometries, good agreement with computational results based on the SC Toolbox, is observed. We also use the reciprocal error estimation method introduced in our earlier paper to validate our numerical results. In particular, exponential convergence, in accordance with the theory of Babu\v ska and Guo, is demonstrated.
\end{abstract}

\maketitle


\section{Introduction}


A bounded Jordan curve in the complex plane divides the extended complex plane $\C_\infty=\C \cup \{\infty\}$ into two domains $D_1$ and $D_2$, whose common boundary it is. One of these domains, say $D_1\,,$ is bounded and the other one is unbounded. The domain $D_1\,$ together with four distinct points $z_1, z_2, z_3, z_4$ in $\partial D_1\,,$ which occur in this order when traversing the boundary in the positive direction, is called a quadrilateral and denoted by $(D_1;z_1, z_2, z_3, z_4) \,$ \cite{Ah,hen,kuh,LV}.

By Riemann's mapping theorem, the domain $D_1\,$ can be mapped conformally onto a rectangle $f\colon D_1  \to (0,1) \times  (0,h) \,$  such that the four distinguished points are mapped onto the vertices of the rectangle $f(z_1)=0\,,$ $f(z_2)=1,$ $f(z_3)=1+ih,$ $f(z_4)=ih\,.$ The unique number $h$ is called the (conformal) modulus of the quadrilateral $(D_1;z_1, z_2, z_3, z_4) \,$ \cite{Ah,hen,kuh,LV}. Apart from its theoretical significance in geometric function theory, the conformal modulus is closely related to certain physical quantities which also occur in engineering applications. In particular, the conformal modulus plays an important role in determining resistance values of integrated circuit networks (see e.g. \cite{ps2,sl}). Similarly, one can map $D_2\,,$ the complementary domain, conformally $g\colon D_2  \to (0,1) \times (0,k) \,$  such that the four boundary points are mapped onto the vertices of the rectangle $g(z_1)=0\,,$ $g(z_2)=1,$ $g(z_3)=1+ik, g(z_4)= ik\,$, reversing the orientation. Again the number $k$ is unique and it is called the exterior modulus of $(D_1;z_1, z_2, z_3, z_4) \,.$ In practice, the computation of both the modulus and the exterior modulus is carried out by using numerical methods such as numerical conformal mapping. Mapping problems involving unbounded domains likewise are related to some well known engineering applications such as determining two dimensional potential flow around a cylinder, or an airfoil.

In the case of domains with polygonal boundary, numerical methods based on the Schwarz-Christoffel
formula have been extensively studied, see \cite{DrTr}. One of the pioneers of numerical conformal mapping was D. Gaier \cite{G}, \cite{P}. The literature and software dealing with numerical conformal mapping
problems is very wide, see e.g. \cite{DrTr} and \cite{ps2}. 
In our earlier paper \cite{hrv} we applied an alternative approach which reduces the problem to the Dirichlet-Neumann problem for the Laplace equation. Thus any software capable of solving this problem may be used. We use the $hp$-FEM method for computing the modulus
of a bounded quadrilateral and here we will apply the same method for the exterior modulus and another method, AFEM \cite{bsv}, for the sake of comparison, as in \cite{hrv}. Our approach also applies to the case of domains bounded
by circular arc boundaries as we will see below. It should be noted that while our method does not require finding the canonical conformal mapping, it is possible to construct the mapping from the potential function. An algorithm, with several numerical examples, is presented in \cite{hqr}.
An alternative to FEM would be to use numerical methods for integral equations.
For recent work on numerical conformal mapping based on such an approach, see Nasser \cite{N}.

In particular, an important example of a quadrilateral $(D_1;z_1, z_2, z_3, z_4) \,$ is the case when $D_1$ is a polygon with $z_1, z_2, z_3, z_4 \,$ as the vertices and its modulus was computed in \cite{hvv} and this formula was also applied in \cite{hrv}. Here we reduce its exterior modulus to the (interior) modulus by carrying out a suitable inversion which keep three vertices invariant and maps the exterior to the interior of a bounded plane region whose boundary consists of four circular arcs.


We apply here three methods to study our basic problem:
\begin{itemize}
\item[(1)] The $hp$-FEM method introduced in \cite{hrv} and its implementation by H.~Hakula.
\item[(2)] The AFEM method of K.~Samuelsson, see e.g. \cite{bsv} and \cite{hrv}.
\item[(3)] The Schwarz-Christoffel Toolbox of T.~Driscoll and N.~Trefethen \cite{Dr,DrTr,t,td}.
\end{itemize}
The methods (1) and (2) are based on a reduction of the exterior modulus problem to the solution of the Dirichlet-Neumann problem for the Laplace equation in the same way as in \cite{Ah} and \cite{hrv} whereas (3) makes use of
numerical conformal mapping methods. Note that \cite{Ah} also provides a connection between the extremal length of a family of curves and its reciprocal, the modulus of a curve family, both of which are widely used in the geometric function theory.

We describe the high-order $p$-, and $hp$-finite element methods and report the results of numerical computation of the exterior moduli of a number of quadrilaterals. In the $p$-method, the unknowns are coefficients of some polynomials that are associated with topological entities of elements, nodes, sides, and the interior.
Thus, in addition to increasing accuracy through refining the mesh, we have an additional refinement parameter, the polynomial degree $p$. For an overview of the $hp$-method, see e.g. Babu\v{s}ka and Suri \cite{bs}. A more detailed exposition of the methods is given in \cite{s,sz}.

Our study is structured according to a few particular cases. We start out with the case when
the quadrilateral is the complement of a rectangle and the vertices are the distinguished
points of the quadrilateral. In this case we have the formula of P.~Duren and J.~Pfaltzgraff
\cite{dp} to which we compare the accuracy of each of the above methods (1)-(3). Then we consider
the problem of minimizing the exterior modulus of a trapezoid with a fixed height $h$ and
fixed lengths for the pair of parallel opposite sides and present a conjecture supported by
our experiments. We also remark that the case of symmetric hexagons can be dealt with the Schwarz-Christoffel transformation and relate its exterior modulus to a symmetry property of the modulus of a curve family. Finally, we study the general case and present comparisons of methods (1)-(3) for this case as well. SC~Toolbox does not have a built in function for computing the exterior modulus. However, we use the function {\tt extermap}, and an auxiliary M\"obius transformation, to map the exterior of a quadrilateral $(D; a,b,c,d)$  conformally onto the upper half-plane so that the boundary points $a,b,c$ and $d$ are mapped to the points $\infty, -1,0$ and $t>0$, respectively. Then the exterior modulus of the
quadrilateral is $\tau(t)/2$, where $\tau$ is the Teichm\"uller modulus function (see \cite{avv} and \ref{mod22} below). We use the MATLAB code from \cite{avv} to compute values of $\tau(t)$, $t>0$.

Our computational workhorse, the $hp$-FEM algorithm implemented in Mathematica, is used in all
cases involving general curved boundaries. We demonstrate that nearly the optimal rate of convergence, $C_1\exp(-C_2 N^{1/3})$ in terms of the number $N$ of unknowns as predicted by the results of Babu\v ska and Guo \cite{bg1}, is attained in a number of tests cases. Our results are competitive with the survey results on $hp$-adaptive algorithms reported by Mitchell and MacClain \cite{Mitchell2}
for the L-shaped domain.

At the end of the paper we present conclusions concerning the performance of these methods
and our discoveries.

\section{Preliminaries}

In this section we give reference results which can be used in obtaining error estimates. We also present some geometric identities which are required in our computations.

\subsection{The hypergeometric function and complete elliptic integrals}
Given complex numbers
$a,b,$ and $c$ with $c\neq 0,-1,-2, \ldots $,
the {\em Gaussian hypergeometric function} is the analytic
continuation to the slit plane $\C \setminus [1,\infty)$ of
the series
\begin{equation} \label{eq:hypdef}
F(a,b;c;z) = {}_2 F_1(a,b;c;z) =
\sum_{n=0}^{\infty} \frac{(a,n)(b,n)}{(c,n)} \frac{z^n}{n!}\,,\:\: |z|<1 \,.
\end{equation}
Here $(a,0)=1$ for $a \neq 0$, and $(a,n)$
is the {\em shifted factorial function}
or the {\em Appell symbol}
$$
(a,n) = a(a+1)(a+2) \cdots (a+n-1)
$$
for $n \in \N \setminus \{0\}$, where
$\N = \{ 0,1,2,\ldots\}$
and the
{\it elliptic integrals} $\K(r),\K'(r)$
of the first kind are
$$
\K(r)=\frac{\pi}{2} F(1/2,1/2;1; r^2),
\qquad
\K'(r)=\K(r'),\text{ and }r'= \sqrt{1-r^2},
$$
and the {\it elliptic integrals} $\E(r),\E'(r)$ of the second  kind are
$$
\E(r)=\frac{\pi}{2} F(1/2,-1/2;1; r^2),
\qquad
\E'(r)=\E(r'),\text{ and }r'= \sqrt{1-r^2}.
$$
Some basic properties of these functions can be found in \cite{avv,olbc}.

\subsection{The modulus of a curve family}
\label{mod22}
For a family of curves $\Gamma$ in the plane, we use the notation $\M(\Gamma)$ for its modulus \cite{LV}.
For instance, if $\Gamma$ is the family of all curves joining the opposite $b$-sides
within the rectangle $[0,a]\times[0,b], a,b>0,$ then $\M(\Gamma)=b/a\,.$ If we consider
the rectangle as a quadrilateral $Q$ with distinguished points $a+ib, ib,0,a$ we also have
$\M(Q;a+ib, ib,0,a)=b/a \,,$ see \cite{Ah,LV} . Given three sets $D,E,F$ we use the notation
$\Delta(E,F; D)$ for the family of all curves joining $E$ with $F$ in $D\,.$

Next consider another example, which is important for the sequel. For $t>0$ let $E= [-1,0]$, $F = [t, \infty)$ and let $\Delta_t$ be the family of curves joining $E$ and $F$ in the upper half-plane ${\mathbb C}_+ = \{ z \in \mathbb{C}: \,\, {\rm Im} \, z > 0  \}\,.$ Then \cite{avv}, we have
\[
\M(\Delta_t)= \tau(t)/2\,;\quad  \tau(t) = 2\frac{ \K(1/\sqrt{1+t})}{\K(\sqrt{t/(1+t)})}\,.
\]

\subsection{The Duren-Pfaltzgraff formula \cite[Theorem 5]{dp}}\label{dpformula} 
For $k \in(0,1)$ write
$$
\psi(k) = \frac{2( \E(k)-(1-k) \K(k)  )}{ \E'(k) - k \K'(k)} \,.
$$
Then $\psi\colon (0,1)\to (0,\infty)$ defines an increasing homeomorphism with limiting
values $0,\infty$ at $0,1$, respectively. In particular, $\psi^{-1}\colon  (0,\infty) \to (0,1)$ is well-defined.
Let $R$ be a rectangle with sides of lengths $a$ and $b$, respectively, and let $\Gamma$
be the family of curves lying outside  $R$ and joining the opposite sides
of length $b\, .$ Then
\begin{equation} \label{dpform}
 \M(\Gamma) = \frac{\K'(k)}{2 \K(k)}\,,\quad \textrm{where}\, \,
k= \psi^{-1}(a/b) \,.
\end{equation}
This formula occurs in different contexts. For instance, W.G. Bickley (\cite{bickley}, (1.17) p. 86) used it in the analysis of electric potentials and
W. von  Koppenfels and F. Stallmann (\cite{VS},  (4.2.31) and (4.2.63))   established it in conformal mapping problems. As far as we know, Duren and
Pfaltzgraff were the first authors to connect this formula with the exterior modulus of a quadrilateral.


\subsection{Mapping unbounded onto bounded domains}
\begin{figure}
\begin{center}
\includegraphics[height=6cm]{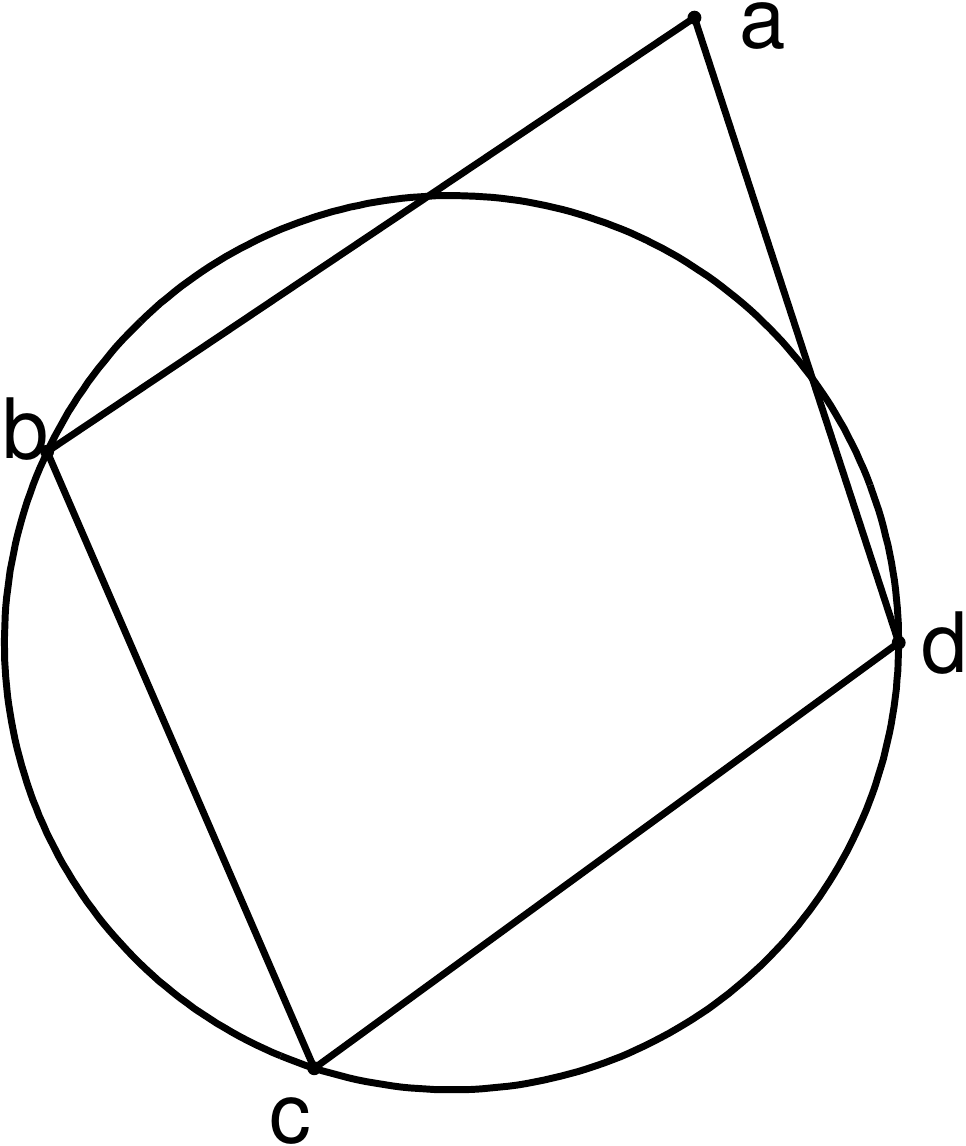}\quad
\includegraphics[height=6cm]{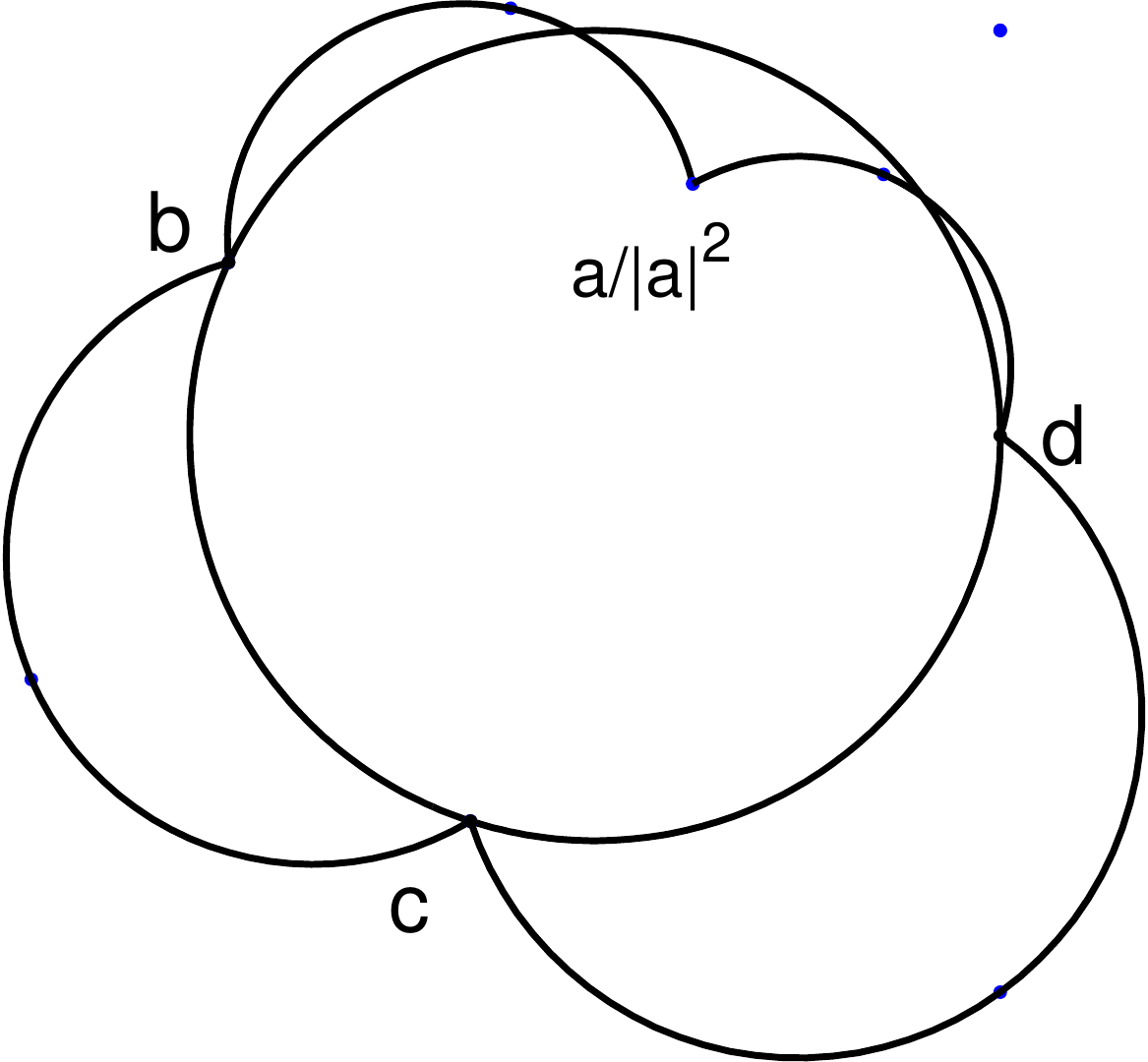}
\end{center}
\caption{Polygonal quadrilateral before (left) and after
(right) the inversion transformation
$z \mapsto z/|z|^2$. Note that the points $b,c,d$ on the unit circle
remain invariant.
}\label{qlat}
\end{figure}

The transformation $z \mapsto z/|z|^2$ maps the complement of the closed unit
disk onto the unit disk. This transformation is an anticonformal mapping
and it maps the complement of a polygonal quadrilateral with vertices $a,b,c,d$ with
$|b|=|c|=|d|=1$ onto a bounded domain, bounded by four circular arcs. Note that
the points $b,c,d$ remain invariant under this transformation. See Figure
\ref{qlat}. Here we also make use of the well-known formula for the center of
the circle through three given points.

\subsection{The Dirichlet-Neumann problem}
The following problem is known as the {\it Dirichlet-Neumann problem}.
Let $D$ be a region in the complex plane whose boundary
$\partial D$ consists of a finite number of regular Jordan
curves, so that at every point, except possibly at finitely many points,
of the boundary a normal is defined. Let $\partial D =A \cup B$
where $A, B$ both are unions of Jordan arcs. Let $\psi_A,\psi_B$ be a real-valued
continuous functions defined on $A,B$, respectively. Find a function $u$
satisfying the following conditions:
\begin{enumerate}
\item
$u$ is continuous and differentiable in
$\overline{D}$.
\item
$u(t) = \psi_A(t),\qquad \textrm{ for all}\,\, t \in A$.
\item
If $\partial/\partial n$ denotes differentiation in
the direction of the exterior normal, then
$$
\frac{\partial}{\partial n} u(t)=\psi_B(t),\qquad \textrm{ for all}\,\, t\in B.
$$
\end{enumerate}

\subsection{Modulus of a quadrilateral and Dirichlet integrals}
One can express the modulus of a quadrilateral $(D; z_1, z_2, z_3, z_4)$
in terms of the solution of the Dirichlet-Neumann problem as follows.
Let $\gamma_j, j=1,2,3,4$ be the arcs of
$\partial D$ between $(z_4, z_1)\,,$ $(z_1, z_2)\,,$ $(z_2, z_3)\,,$
$(z_3, z_4),$ respectively. If $u$ is the (unique) harmonic solution of
the Dirichlet-Neumann problem with boundary values of $u$ equal to $0$ on
$\gamma_2$, equal to $1$ on $\gamma_4$ and with $\partial u/\partial n =
0$ on $\gamma_1 \cup \gamma_3\,,$ then by \cite[p. 65/Thm 4.5]{Ah}:
\begin{equation}
\label{qmod}
\M(D;z_1,z_2,z_3,z_4)=
\iint_D |\triangledown
u|^2\,dx\,dy.
\end{equation}

\subsection{The reciprocal identity}
Given a quadrilateral $Q=(D; z_1, z_2, z_3, z_4)$ we call sometimes
$\tilde{Q}=(D;  z_2, z_3, z_4,z_1)$ the conjugate quadrilateral. It
a simple basic fact that
\begin{equation} \label{recipidty}
\M(Q) \M(\tilde{Q}) =1\,.
\end{equation}
It was suggested in \cite{hvv} and \cite{hrv} that the quantity
\begin{equation} \label{recipidty2}
r(Q) =|\M(Q)\M(\tilde{Q}) -1|\,
\end{equation}
might serve as a useful error characteristic. We will continue to
use this also in our work.

\subsection{The $hp$-FEM method and meshing} \label{meshing}

In this paper, we use the $hp$-FEM method for computing for the exterior modulus of a quadrilateral. 
For a general description of our method, see \cite{hrv}. 
Proper treatment of corner singularities is handled with the following two-phase algorithm,
typically recursive,
where triangles can be replaced by quadrilaterals or a mixture of both:
\begin{enumerate}
\item Generate an initial mesh (triangulation) where the corners are
isolated with a fixed number of triangles depending on the interior angle, $\theta$ so that the refinements can be carried out independently:
\begin{enumerate}
\item $\theta \leq \pi/2$: one triangle,
\item $\pi/2 < \theta \leq \pi$: two triangles, and
\item $\pi < \theta$: three triangles.
\end{enumerate}
\item Every triangle attached to a corner is replaced by a refinement,
where the edges incident to the corner are split as specified by the scaling factor
$\alpha$. This process is repeated recursively until the desired nesting level $\nu$
is reached. The resulting meshes are referred to as $(\alpha,\nu)$-meshes. 
Note that the mesh may include quadrilaterals after refinement.
\end{enumerate}
Since the choice of the initial mesh affects strongly the refinement
process, it is advisable to test with different choices. Naturally,
one would want the initial mesh to be minimal, that is, having the
smallest possible number of elements yet providing support for the
refinement. This is why initial meshes are sometimes referred to as minimal
meshes.

\begin{figure}
\begin{center}
\includegraphics[width=.45\textwidth]{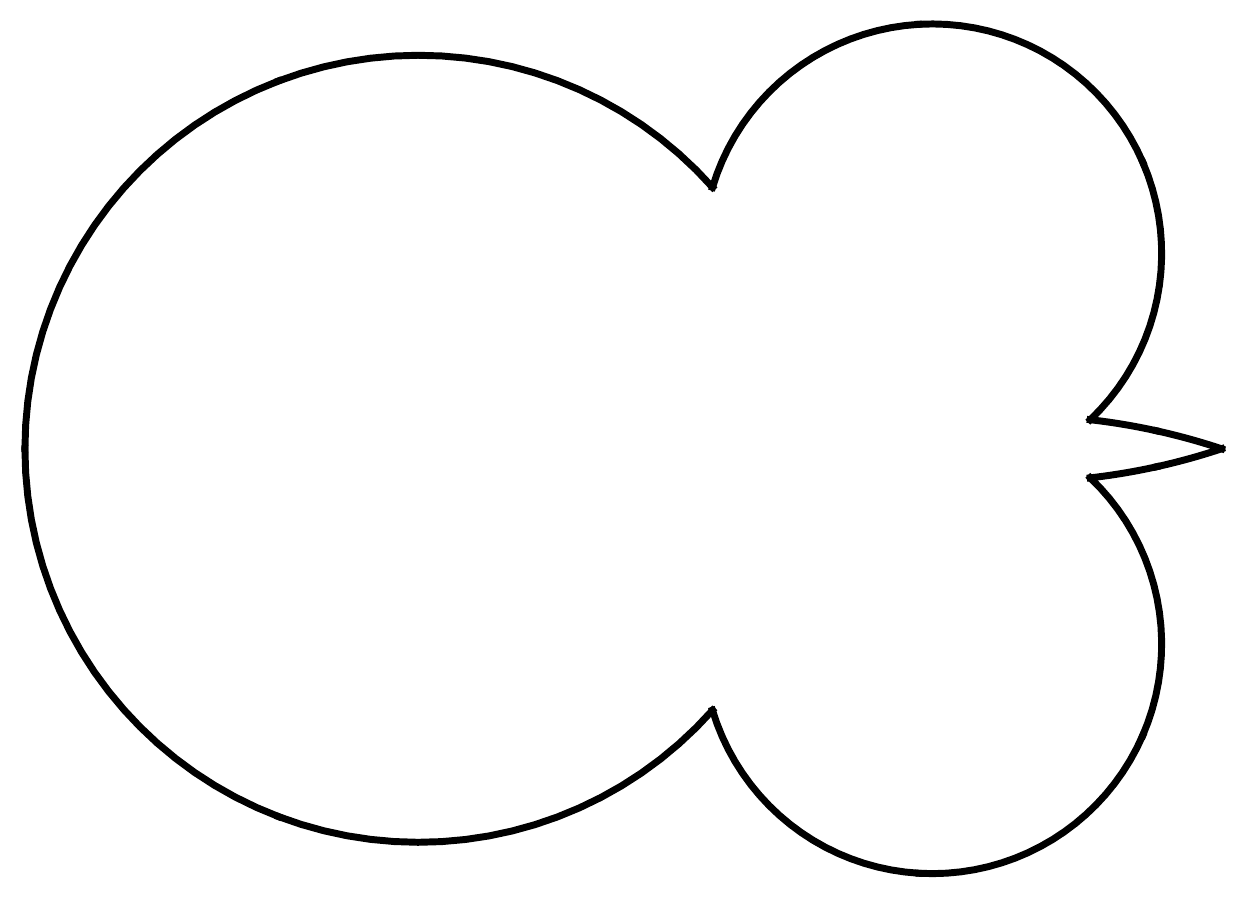}
\includegraphics[width=.45\textwidth]{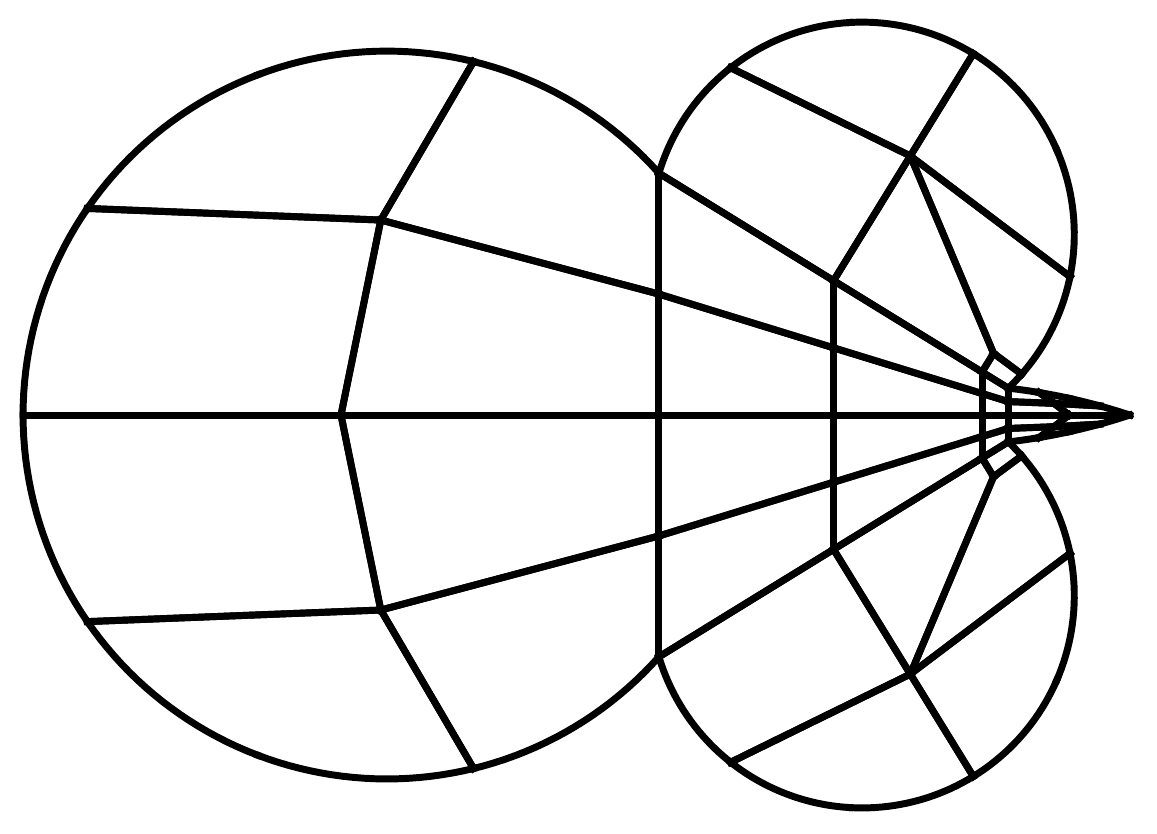}
\end{center}
\caption{A sample geometry and the corresponding initial mesh. Note
  that the three-element -rule is satisfied at every corner.}\label{mesh:initialconfiguration}
\end{figure}
\begin{figure}
\begin{center}
\includegraphics[width=.3\textwidth]{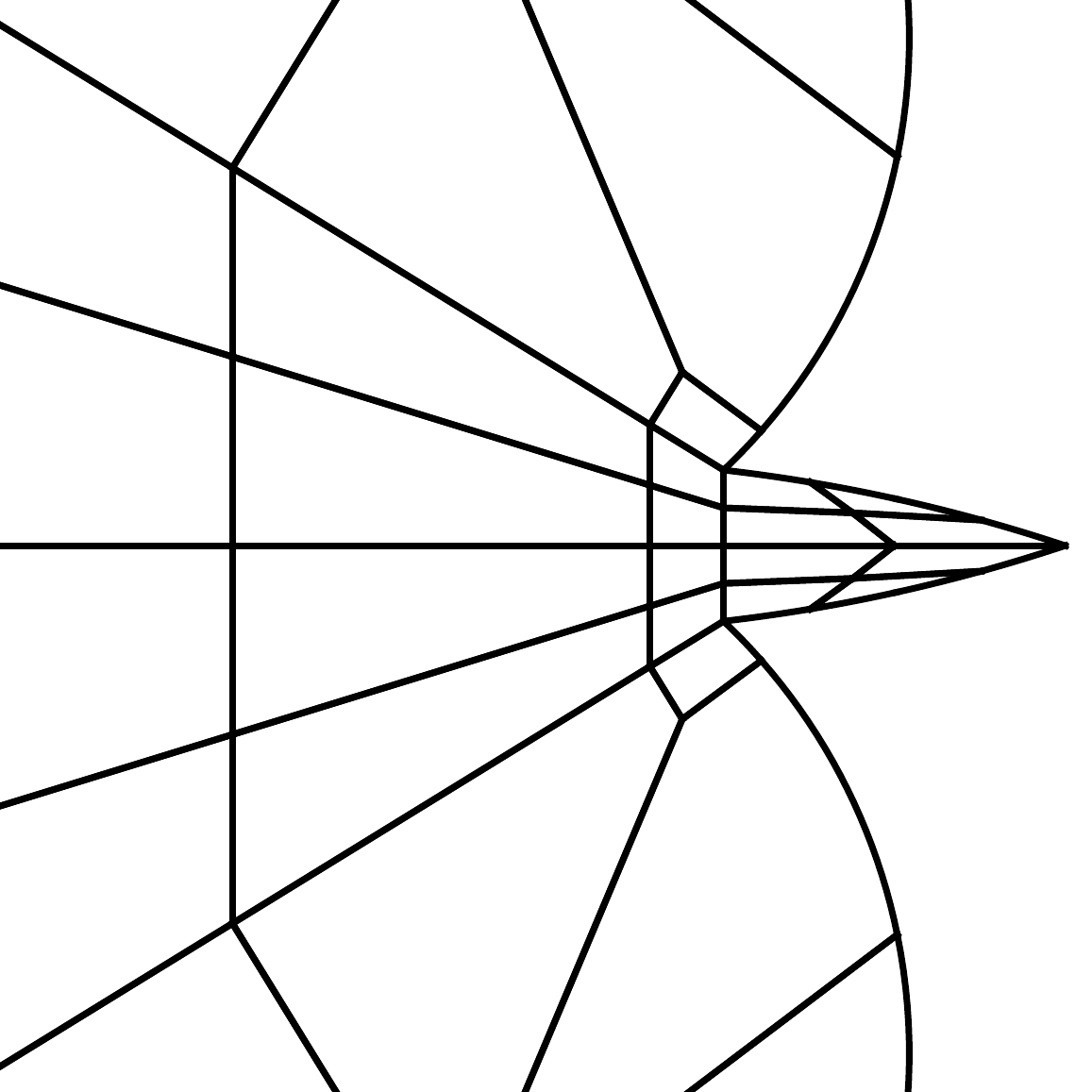}
\includegraphics[width=.45\textwidth]{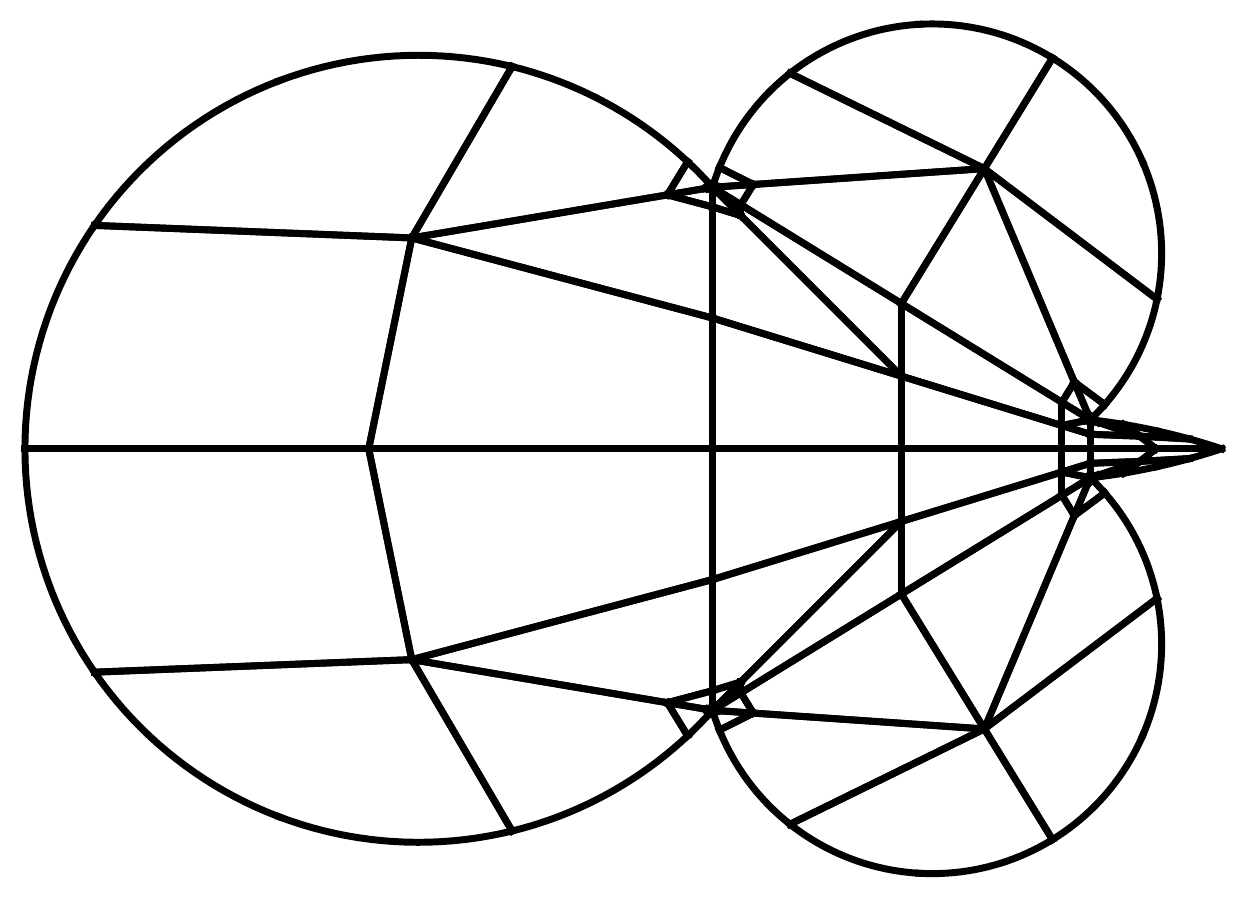}
\end{center}
\caption{A detail of the initial mesh and the final $(0.15, 14)$-mesh
  used in the actual computation.}\label{mesh:detail}
\end{figure}
In Figure~\ref{mesh:initialconfiguration} a challenging example is
shown. In this case the large variation of the edge lengths is
addressed by adding a refinement step to the construction of the
initial mesh. A detail of the initial mesh is given in
Figure~\ref{mesh:detail} along with the final mesh.

\section{The case of a rectangle}

The first tests with the $hp$-FEM software were made for the case of the exterior
modulus of a rectangle and checked against the Duren-Pfaltzgraff formula (\ref{dpform}).
For a convenient parametrization of the computation, the vertices of the rectangle
were chosen to be the points $1, e^{it}, -1, -e^{it}, t \in (0, \pi/2]$ of the unit
circle. In this case, the "interior" modulus of the rectangle is $\tan (t/2)\,.$
It is equal to the modulus of the family of curves joining the sides $[1, e^{it}]$
and $[-1, -e^{it}]$ and lying in the interior of the rectangle. The formula (\ref{dpform})
now gives the corresponding exterior modulus as
$$
\frac{\K'(k)}{2 \K(k)} \,,\qquad k= \psi^{-1}\bigg(\frac{1}{\tan(t/2)}\bigg).
$$

\begin{figure}
\centering
\subfloat[$k=1$]{\includegraphics[width=0.3\textwidth]{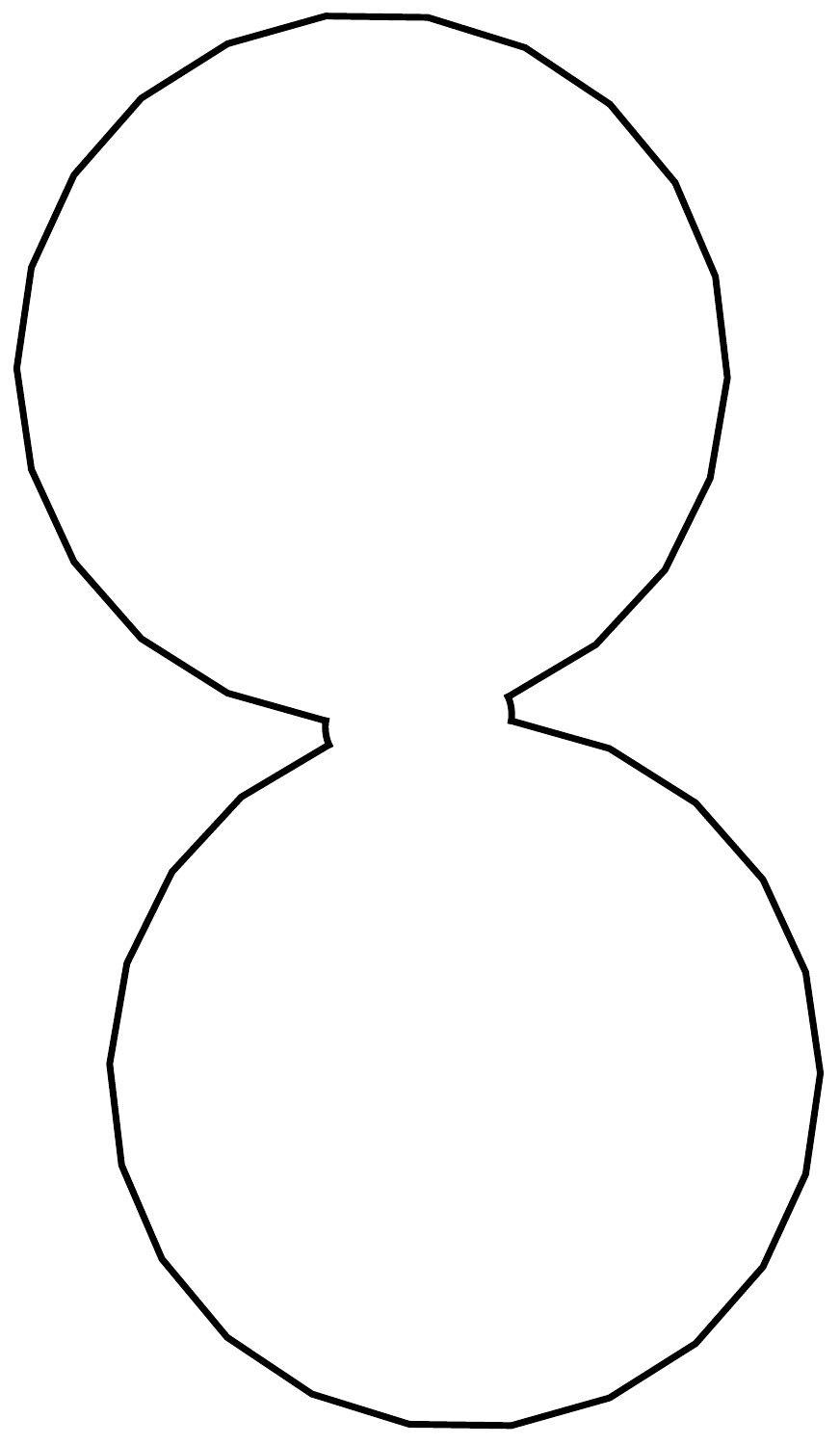}}\quad
\subfloat[$k=2$]{\includegraphics[width=0.3\textwidth]{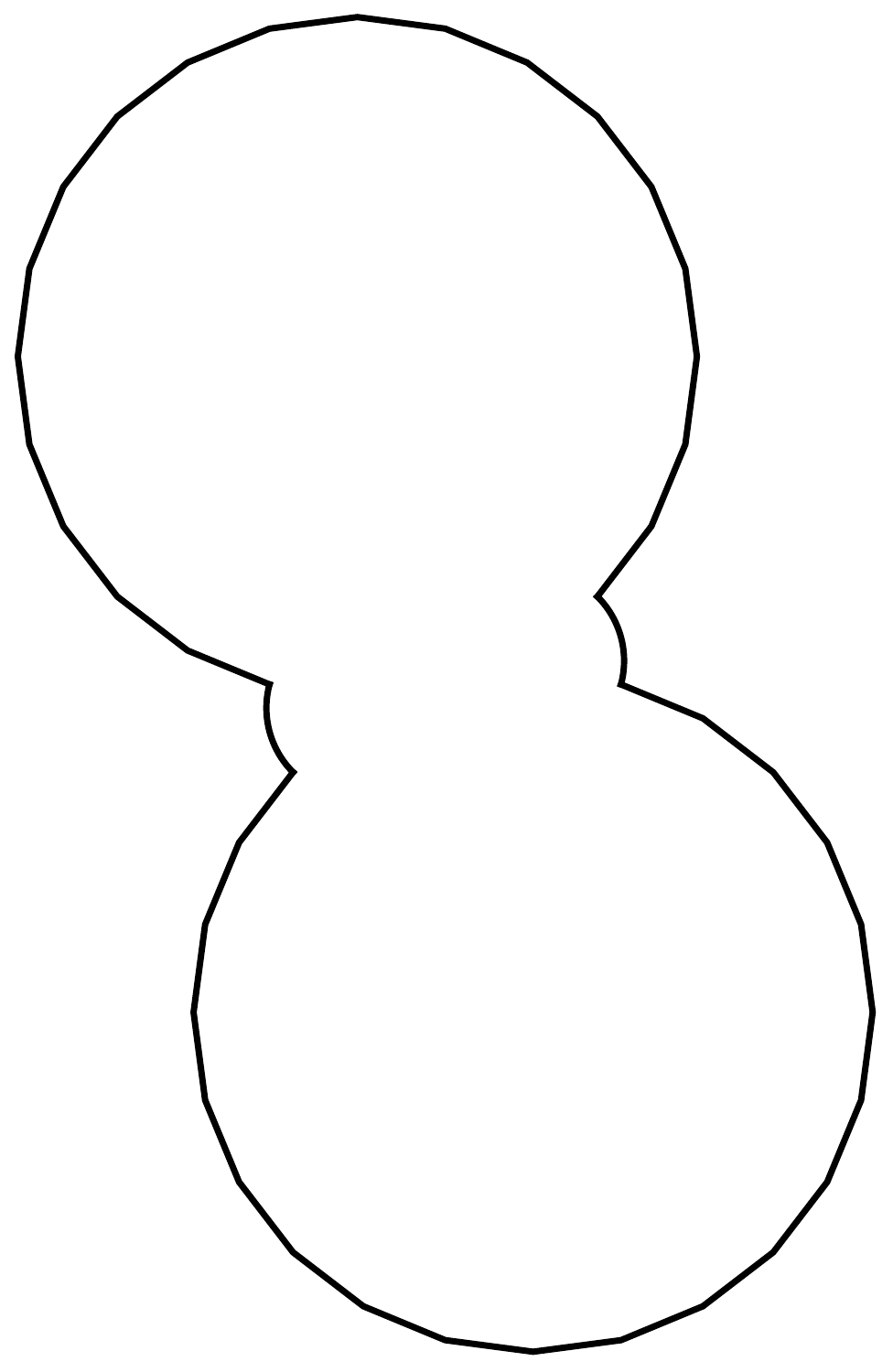}}\quad
\subfloat[$k=3$]{\includegraphics[width=0.3\textwidth]{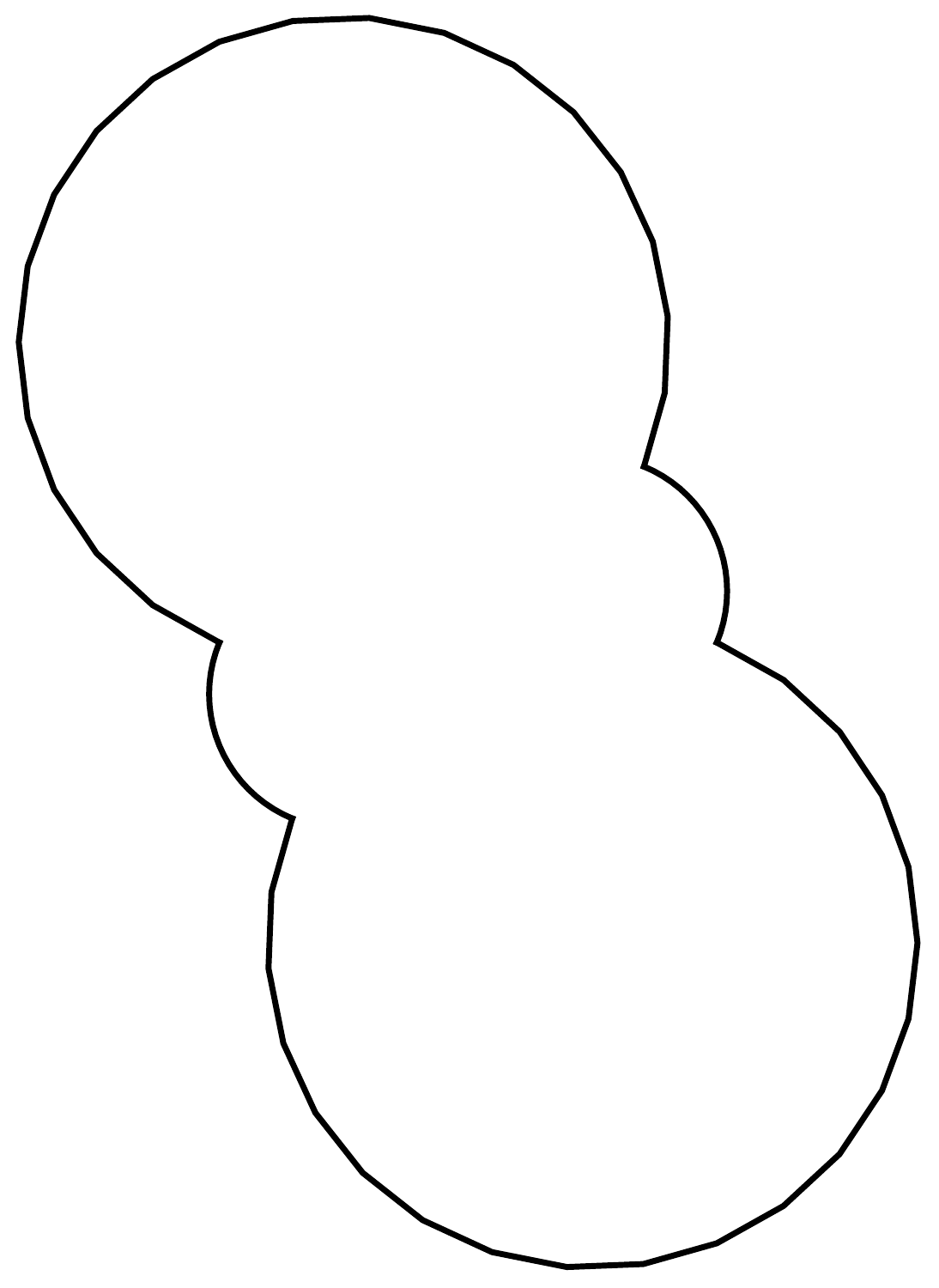}}\\
\subfloat[$k=4$]{\includegraphics[width=0.3\textwidth]{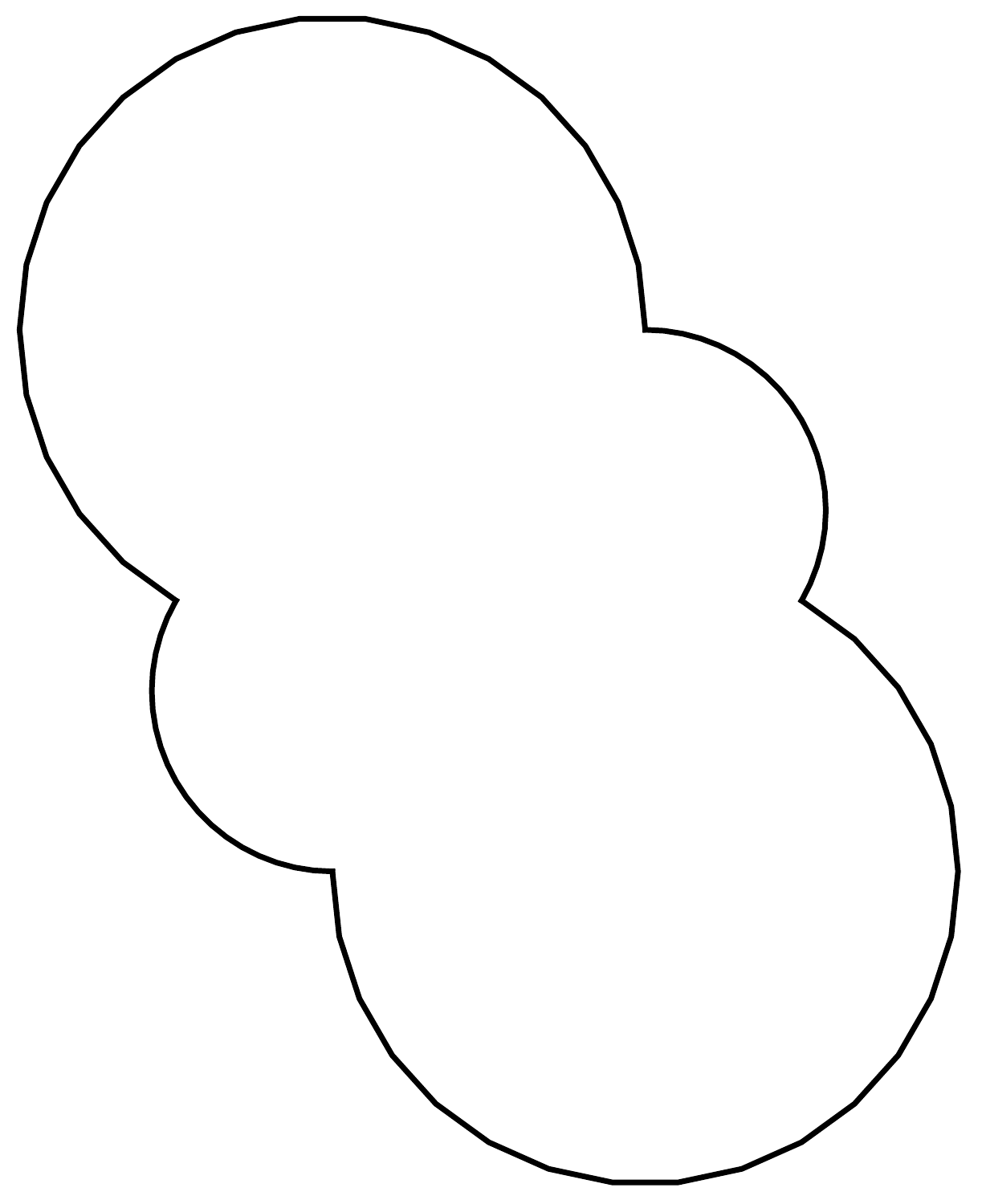}}\quad
\subfloat[$k=5$]{\includegraphics[width=0.3\textwidth]{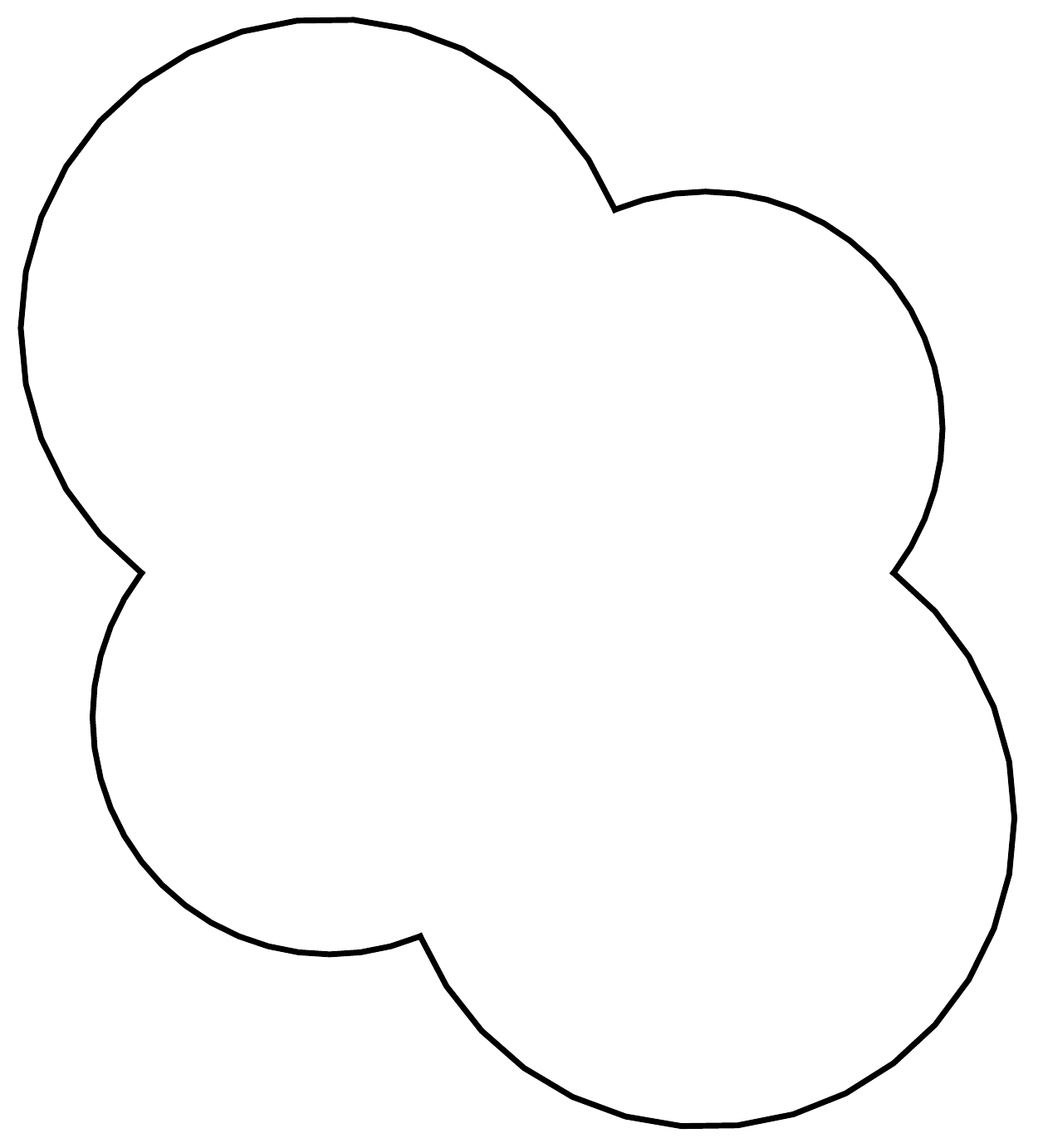}}\quad
\subfloat[$k=6$]{\includegraphics[width=0.3\textwidth]{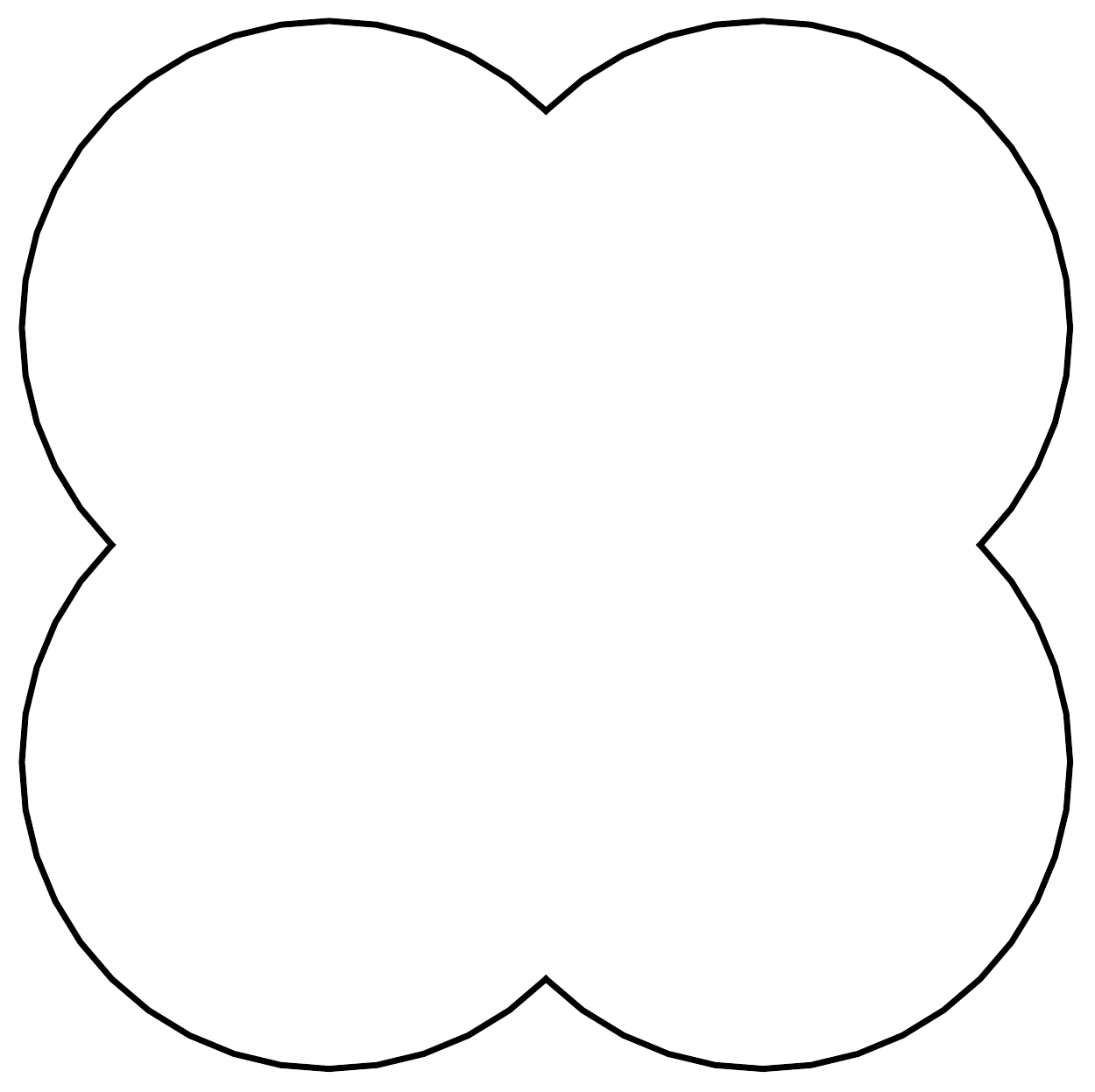}}
\caption{Circular arc domains used for the $hp$-FEM computations
of the values in Table \ref{table101217}. The scale varies from picture to picture.
}\label{domain101217}
\end{figure}

For the computation, we carried out the inversion $z \mapsto 1/\overline{z}$ in the
unit circle which keep all the points of the unit circle fixed and transforms the
exterior modulus problem for the rectangle to the "interior" modulus problem of
a plane domain bounded by four circular arcs, see Figure \ref{domain101217}. These circular arcs are the images
of the sides of the rectangle under the inversion. The results turned out to be
quite accurate, with a typical relative error of the order $10^{-10}\,.$

\begin{table}[ht]
\caption{Exact values of the moduli of $Q(1, e^{it}, -1, -e^{it})$ given by (\protect{\ref{dpform}}) and the errors of computational results
of the $hp$-method, $p=20,$ the AFEM method and the SC Toolbox. The errors are
obtained by comparing with the exact formula
(\ref{dpform}). The errors are given as {$\lceil\log_{10}|\mathrm{error}|\rceil$}.} \label{table101217}  
\begin{tabular}{|c|c|c|c|c|}
\hline
$k$  & exact($t=k \pi/12$) & Error[hpFEM]  & Error[AFEM]& Error[SCT]\\
\hline
$1$ &  $1.50290233467$ & $-9$ &  $    -6$ & $ -9$\\
$2$ &  $1.31044063554$ & $-9$ &  $  -6$ & $ -9$\\
$3$ &  $1.20035166917 $& $-9$ &  $    -6$  & $ -10$\\
$4$ &  $1.12114255114$ & $-10$ &  $ -6$ & $-9 $\\
$5$ & $1.05681535228$ & $-10 $ &  $ -6$ & $ -13$\\
$6$ &  $1. $ &  $-10         $ &  $ -6$ & $ -15$\\
\hline
\end{tabular}
\end{table}

 \begin{figure}
\centering
\includegraphics[width=0.35\textwidth]{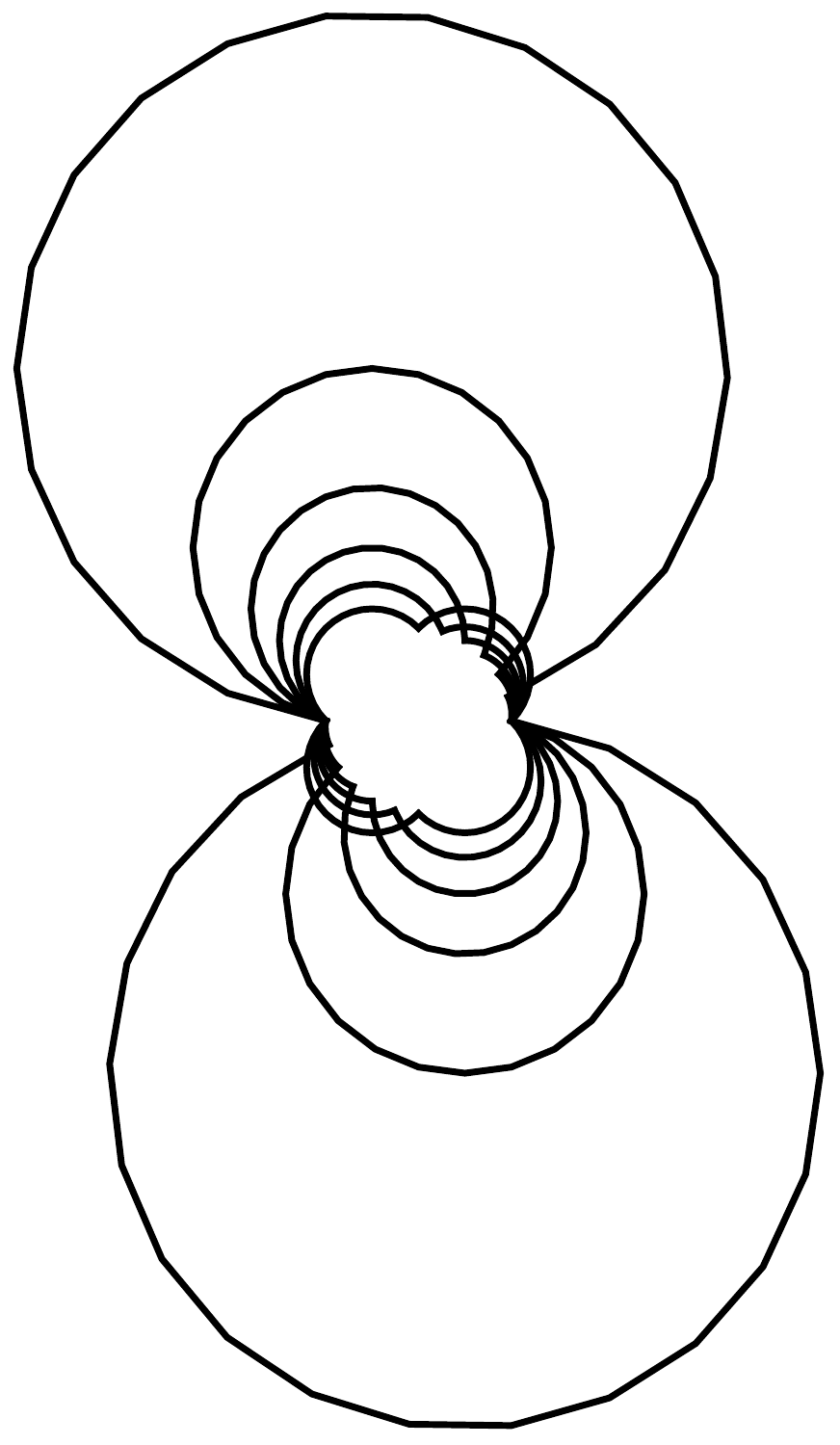}
\caption{The circular arc domains of Figure \ref{domain101217} in the same scale.
}\label{fig:quadsequencecombination}
\end{figure}

\section{Side sliding conjecture}

\subsection{The side sliding problem}\label{myslide}
Consider the problem of finding the minimal exterior modulus of the polygonal
quadrilateral with vertices $0,1,a=t+ ih,b=t-s+ih$ when $h ,s >0$ are fixed
and $t$ varies. We consider the question of computing the modulus of the family $\Gamma$ of curves joining the opposite sides $[1,a]$ and $[b,0]$ outside the quadrilateral.
Our first step is to reduce the problem to an equivalent problem such that three of
the points are on the unit circle.
Note that this setting is valid only if $z_0$ is inside the
quadrilateral. Indeed, for every choice of $h$ and $s$ this condition
defines an upper limit for the value of $t$.

\subsection{Side sliding conjecture}\label{slideconj} 
 The least valueof the exterior modulus is attained when  $t=(1+s)/2\,.$ For  $t\le (1+s)/2\,$
the modulus is a decreasing function of $t\,.$

\subsection{Numerical experiments on side sliding conjecture} \label{numslideconj}
\begin{figure}
\centering
\includegraphics[width=0.55\textwidth]{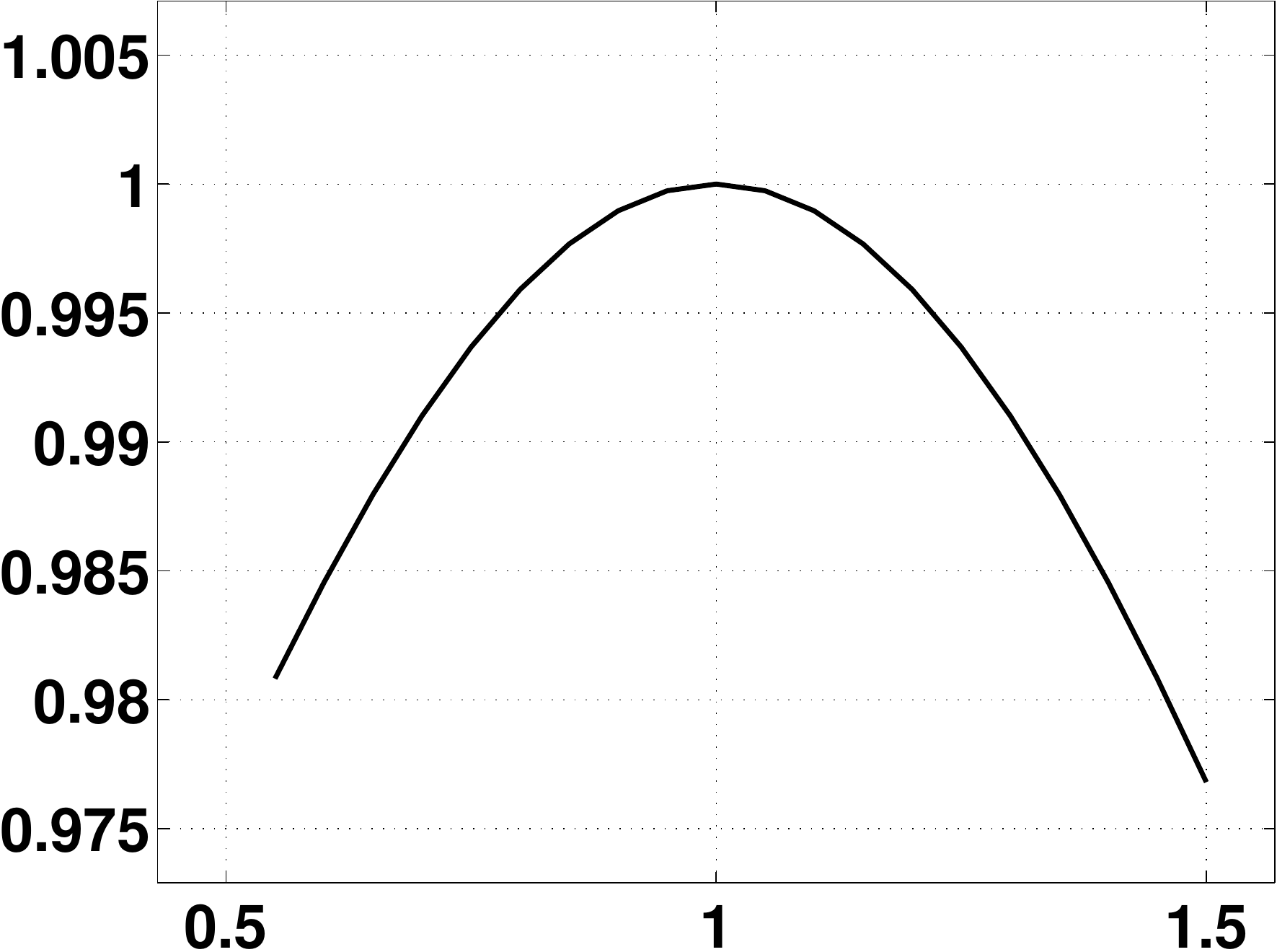}
\caption{Side Sliding Conjecture: Dependence of the exterior modulus
on parameter $t$ with $h=1, s=1$. Maximum is reached at $t=(1+s)/2=1$, as predicted by the conjecture. }\label{fig:sidesliding}
\end{figure}

In Figure~\ref{fig:sidesliding} we show a graph of the exterior
module as a function of the parameter $t\in [0.5,2.5]$, when
$h=1, s=2$.
The computation was carried out with SC Toolbox, $hp$-FEM, and AFEM and for the
range of computed values, the respective graphs coincide.
For the SC Toolbox and the $hp$-FEM the reciprocal estimate for the error was smaller than
$10^{-8}$ and for AFEM $10^{-5}$.

\section{The case of a symmetric hexagon}

Suppose that $Q(a,b,0,1)$ is a quadrilateral in the upper half plane. Then the closed polygonal
line $a,b,0,\overline{b},\overline{a},1,a$ defines a hexagon $H= Q \cup \overline{Q}$ symmetric with respect to the
real axis. Map the complement of $H$ onto $\mathbb{C}\setminus \{-1-t, 1+t\}$ by a conformal
map $h$ such that $h(0)=-1-t, h(b)= h(\overline{b})=-1, h(a)=h(\overline{a})=1, h(1)= 1+t$
where $t>0$ depends on the point configuration $a,b,0,1\,.$  It is clear by symmetry that
\begin{equation}
2 \M(\Delta^+) = \M(\Delta)
\end{equation}
where
$$\Delta=\Delta\big([-1-t,-1],[1,1+t]; \mathbb{C}\big)\, \, \textrm{and} \,\,
\Delta^+=\Delta\big([-1-t,-1],[1,1+t]; \{  z:  \im z >0 \} \big) \,.$$
Because of the conformal invariance
of the modulus we also have
\begin{equation} \label{qlatdupl}
2 \M\big(h^{-1}(\Delta^+)\big) = \M\big(h^{-1}(\Delta)\big) \,.
\end{equation}

Applying this formula with \eqref{dpform} we see that
\begin{equation} \label{dupl}
\M(\Gamma_+)= \frac{\K'(k)}{4\, \K(k)}\,, \quad k = \psi^{-1}\bigg(\frac{1}{2h}\bigg)\,,
\end{equation}
where for $h>0$,
\begin{equation}
\Gamma_+= \Delta\big([0,ih], [1, 1+ih]; {\mathbb C}_+ \setminus [0,1]\times [0, h]\big) \,.
\end{equation}

This formula can be checked by using the SC Toolbox to construct the above conformal
mapping $h\,.$ The tests we carried out for $h=0.2,0.3,0,4$ and $0.5$ In these cases the reciprocal estimate for error was  smaller than $10^{-9}\,.$

\section{General quadrilateral}

The exterior modulus of the quadrilateral $Q$ with vertices $a,b,c,d$
is considered in this section, i.e., we compute
$$
\iint_Q \, |\nabla u|^2\, dx\,dy
$$
over the complement of the quadrilateral when $u$ is the solution of the Laplace
equation in the complement of the quadrilateral with Dirichlet values $1$ and $0$
on the sides $[b,c]$ and
$[d,a]\,,$ respectively, and the Neumann value $0$ on the sides $[a,b]$ and $[c,d]\,.$
Here we allow the boundary of the quadrilateral $\partial Q$, be a parametrized curve
$\gamma(t)$, $t\in [-1,1]$.

\begin{figure}
\centering
\subfloat[Exterior domain with radius $> 1000000$.]{\includegraphics[width=0.4\textwidth]{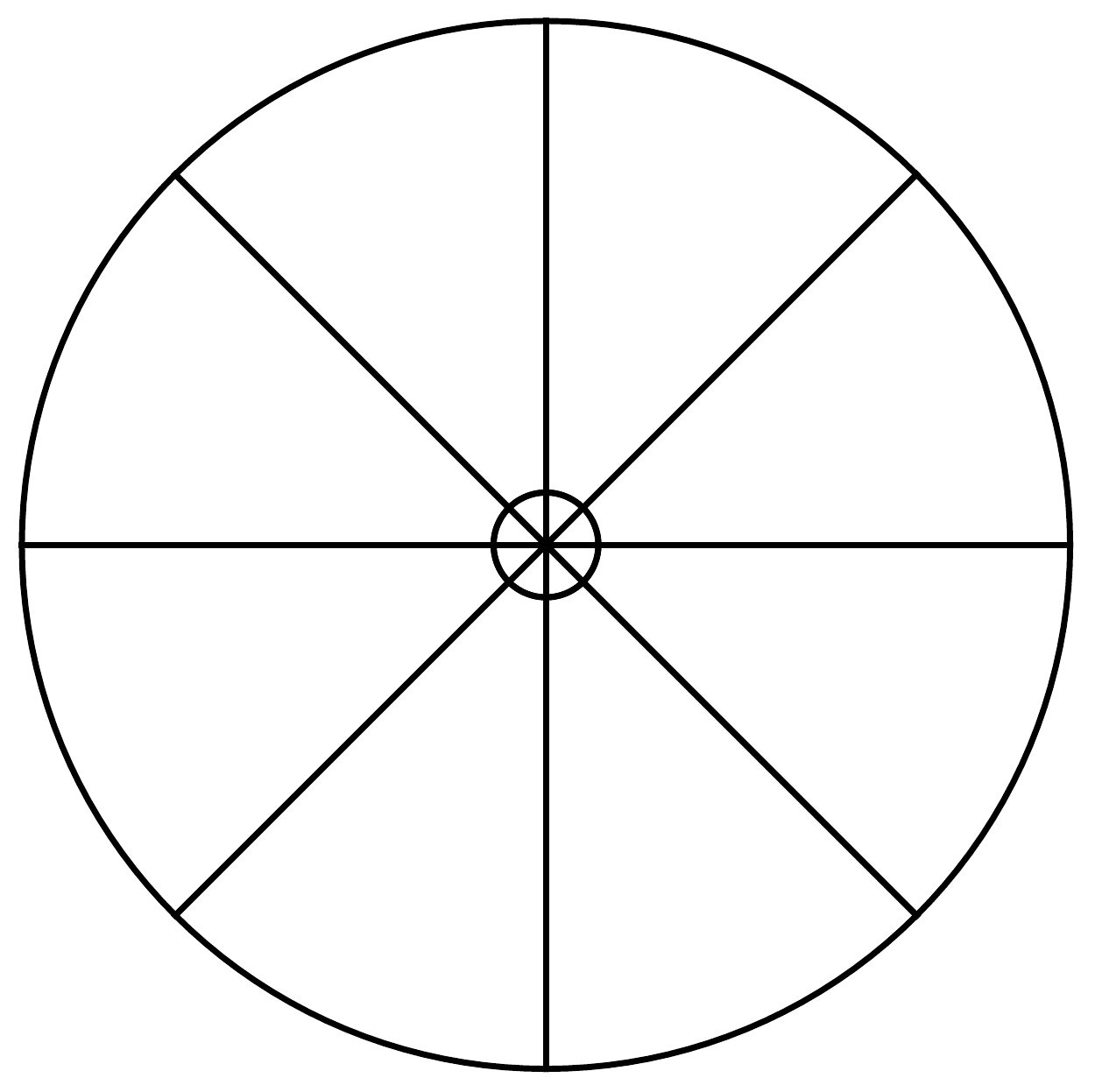}}\\
\subfloat[Zoom of the mesh in the case of a circle.]{\includegraphics[width=0.4\textwidth]{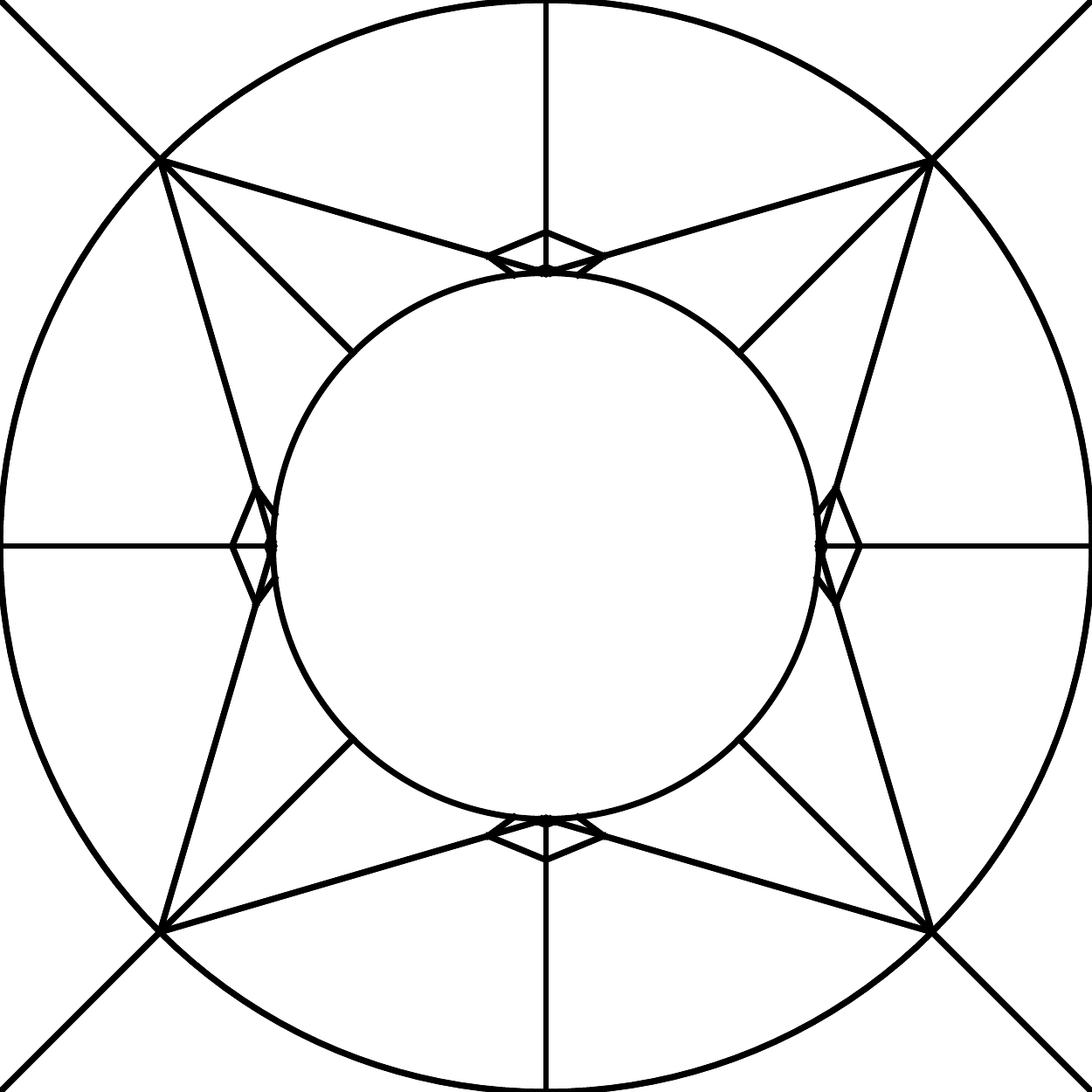}}\qquad
\subfloat[Zoom of the mesh in the case of a square.]{\includegraphics[width=0.4\textwidth]{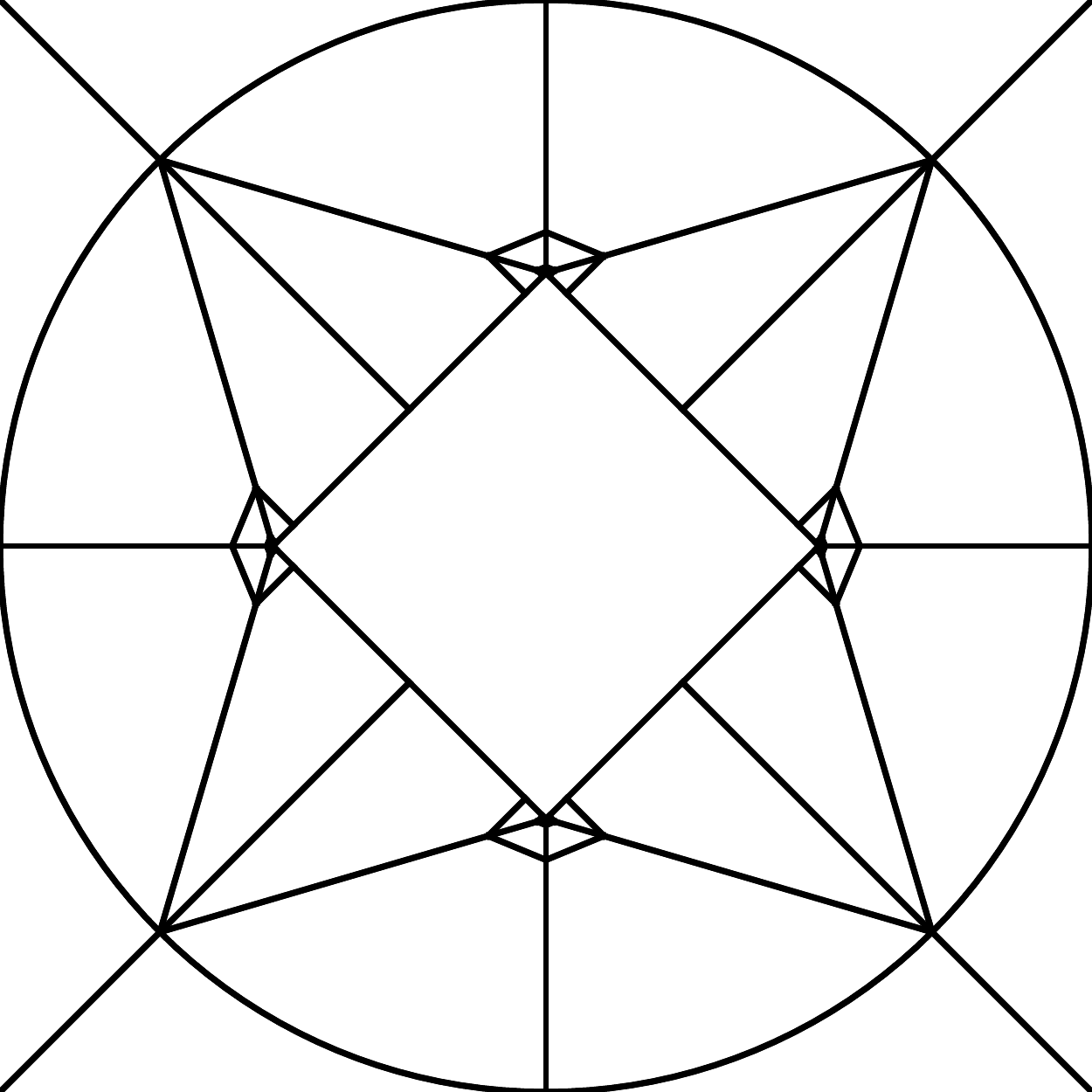}}\\
\subfloat[Zoom of the potential in the case of a circle. Reciprocal
error $\sim 8.8\cdot10^{-10}$]{\includegraphics[width=0.4\textwidth]{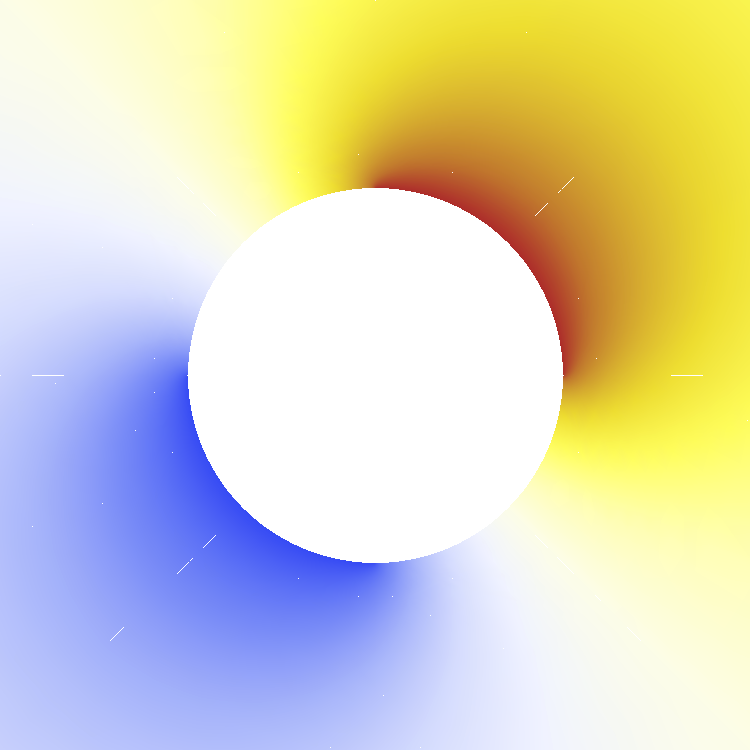}}\qquad
\subfloat[Zoom of the potential in the case of a square. Reciprocal
error $\sim 6.3\cdot10^{-10}$]{\includegraphics[width=0.4\textwidth]{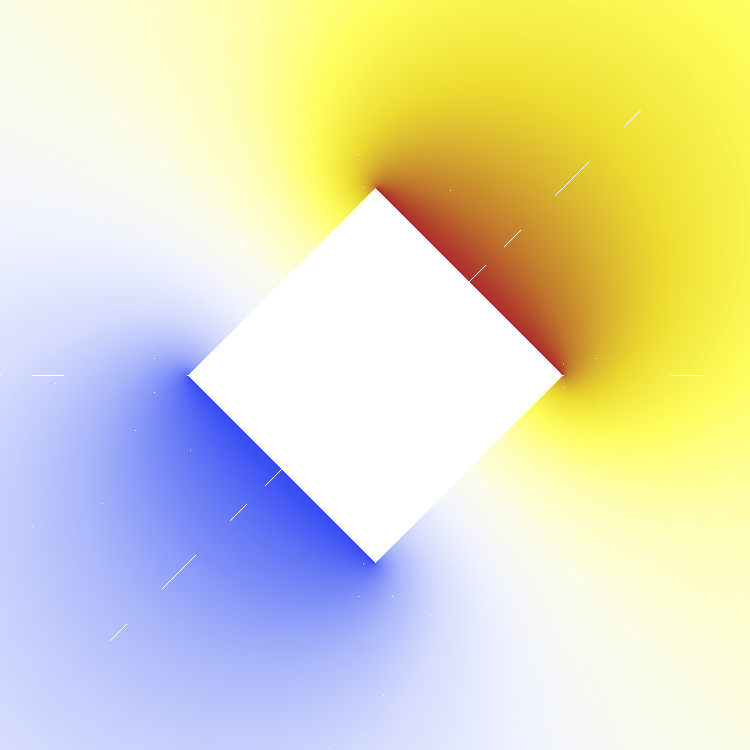}}\\
\caption{Exterior modulus over the exterior domain.}\label{fig:extdom}
\end{figure}

In Figure~\ref{fig:extdom} an overview of the standard FEM approach is given.
Using higher-order elements one can stretch the domain without introducing significant
number of elements. Singularities at the corner point must be accounted for in the grading
of the mesh.
Since both the circle and the square cases are symmetric, the exterior modulus is
exactly 1, and furthermore the potential value at infinity or the far-field value
is exactly 1/2.

\subsection{Quadrilaterals A and B}
In Tables \ref{tbl:extmodcomparisonA}, \ref{tbl:extmodcomparisonB}, and
\ref{tbl:infinitycomparison} results on two polygonal quadrilaterals
\begin{itemize}
\item Quadrilateral A: $\{0, 1, ({28}/{25}, {69}/{50}), (-{19}/{25} , {21}/{25})\}$,
\item Quadrilateral B: $\{0, 1, ({42}/{25}, 4), (-{3}/{25} , {21}/{25})\}$,
\end{itemize}
are presented. The exterior modulus has been computed using three methods as
an equivalent interior modulus problem, and also in truncated domain.
In the interior case, both SC Toolbox and $hp$-FEM give similar results, but AFEM
in its standard setting does not reach the desired accuracy. This is probably due to
the adaptive scheme failing in the presence of cusps in the domain. Tables \ref{tbl:extmodcomparisonA}
and \ref{tbl:extmodcomparisonB} indicate that large exterior angles are the most
significant source of errors in the FEM solutions, as expected. In the rather benign
setting of the Quadrilateral A, SC Toolbox and both the internal and external $hp$-FEM
versions have the same accuracy, but in the case of Quadrilateral B, we see gradual loss
of accuracy in the FEM solutions.

\subsection{Quadrilaterals C and D}
Finally, we consider two flower domains, that is, quadrilateral domains with the boundary
$\gamma(t) = r(t) e^{it}\,, \, r(t) = 4/5 + (1/5) \cos(n \pi t)$ and corners at $t=-1,-1/4,0,1/2$.
For the Quadrilateral C we choose $n=4$ and for D we choose $n=8$.
These domains have the useful property that the exterior problem can easily be
converted to a corresponding interior problem of the domain with boundary $1 / \overline{\gamma(t)}$.
Since these domains cannot be handled using the SC Toolbox, we take the interior
solution as the reference. Tables \ref{tbl:extmodcomparisonC}
and \ref{tbl:extmodcomparisonD} show that we can obtain results of high accuracy also
in traditionally challenging domains.

\begin{figure}
\centering
\subfloat[Contours of the inner problem. Origin is indicated with a dot.]{\includegraphics[width=0.4\textwidth]{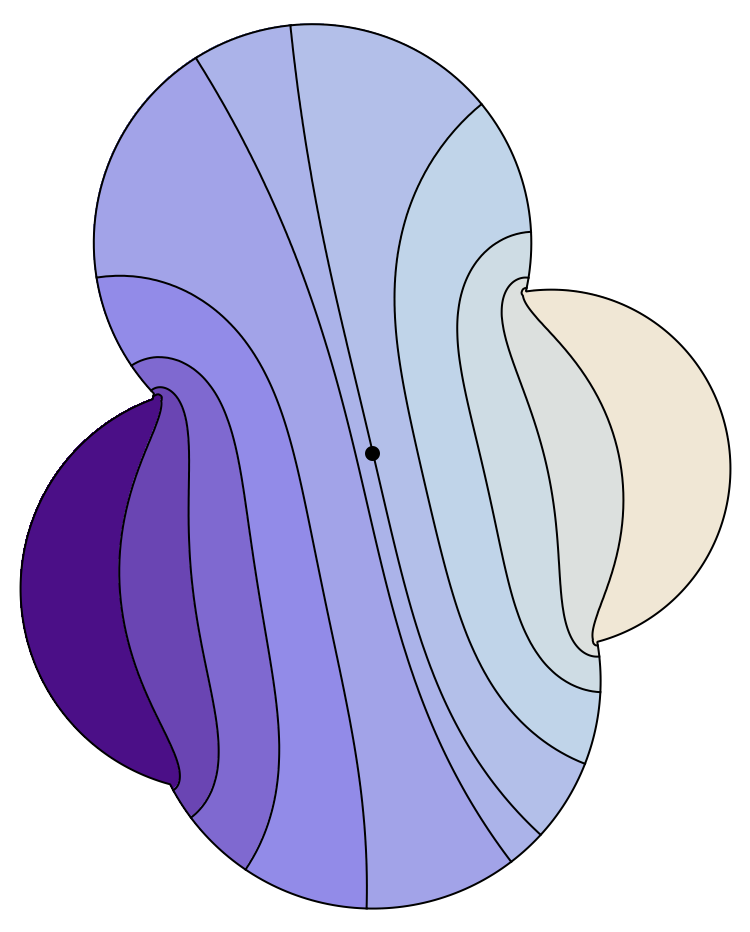}}
\subfloat[Contours of the outer problem. Note the contours extending to infinity.]{\includegraphics[width=0.4\textwidth]{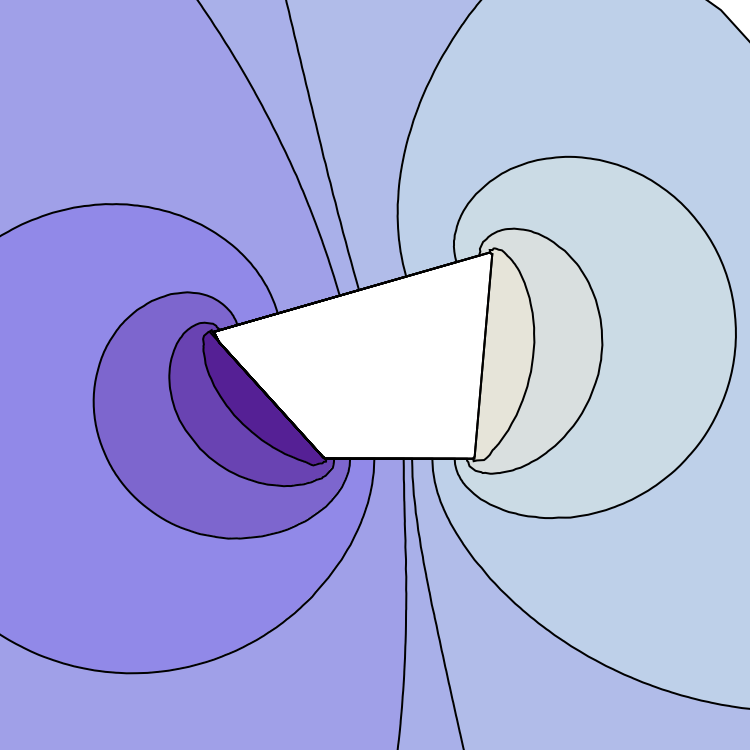}}
\caption{Quadrilateral A: Correspondence of the potential contours between the inner (A) and
outer (B) solutions. Shown are the potential levels $u(z) = 0,1/10,\ldots,1$, and $u(0)$. Corresponding contours can be identified by matching the shadings of
the regions in between.}\label{fig:inneroutercomparisonA}
\end{figure}

\begin{figure}
\centering
\subfloat[Contours of the inner problem. Origin is indicated with a dot.]{\includegraphics[width=0.5\textwidth]{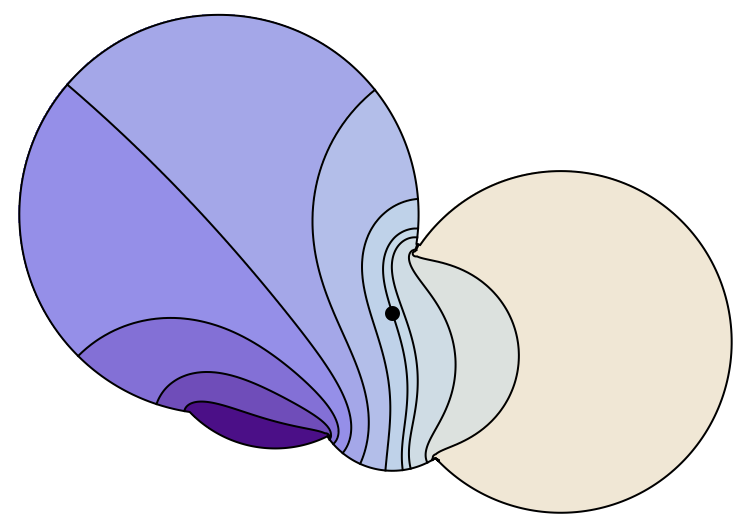}}
\subfloat[Contours of the outer problem. Note the contours extending to infinity.]{\includegraphics[width=0.35\textwidth]{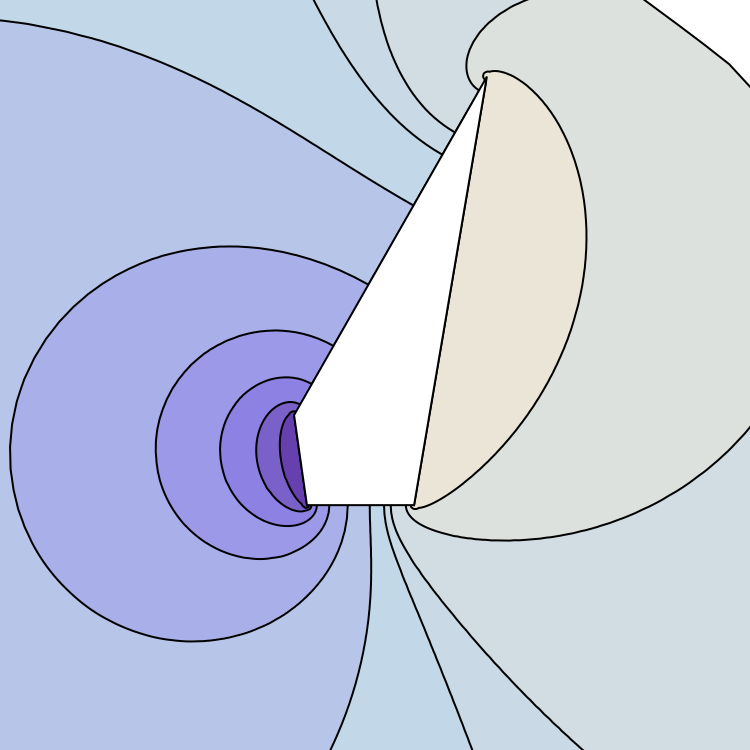}}
\caption{Quadrilateral B: Correspondence of the potential contours between the inner (A) and
outer (B) solutions. Shown are the potential levels $u(z) = 0,1/10,\ldots,1$, and $u(0)$. Corresponding contours can be identified by matching the shadings of
the regions in between.}\label{fig:inneroutercomparisonB}
\end{figure}

\begin{figure}
\centering
\subfloat[Potential field of the inner problem.]{\includegraphics[width=0.4\textwidth]{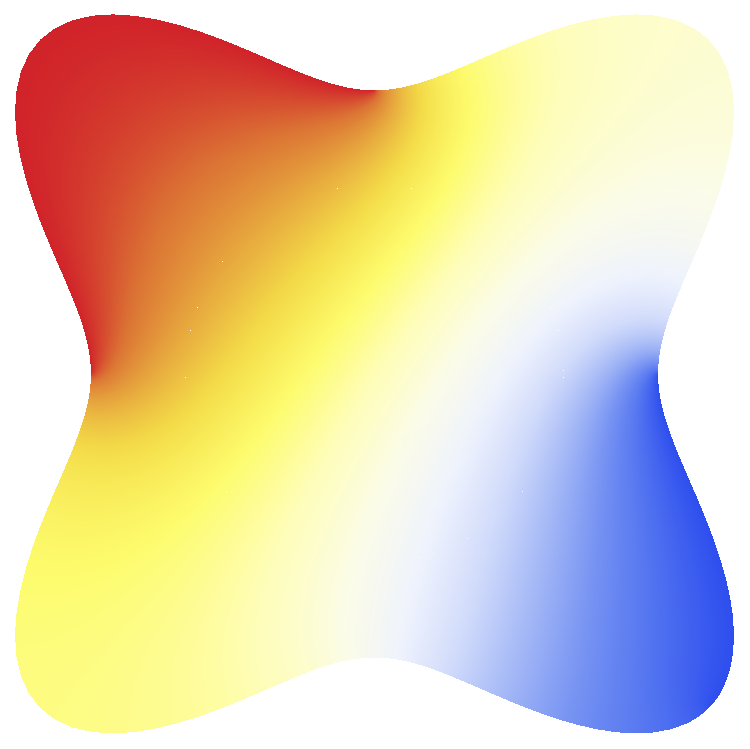}}
\subfloat[Potential field of the outer problem.]{\includegraphics[width=0.4\textwidth]{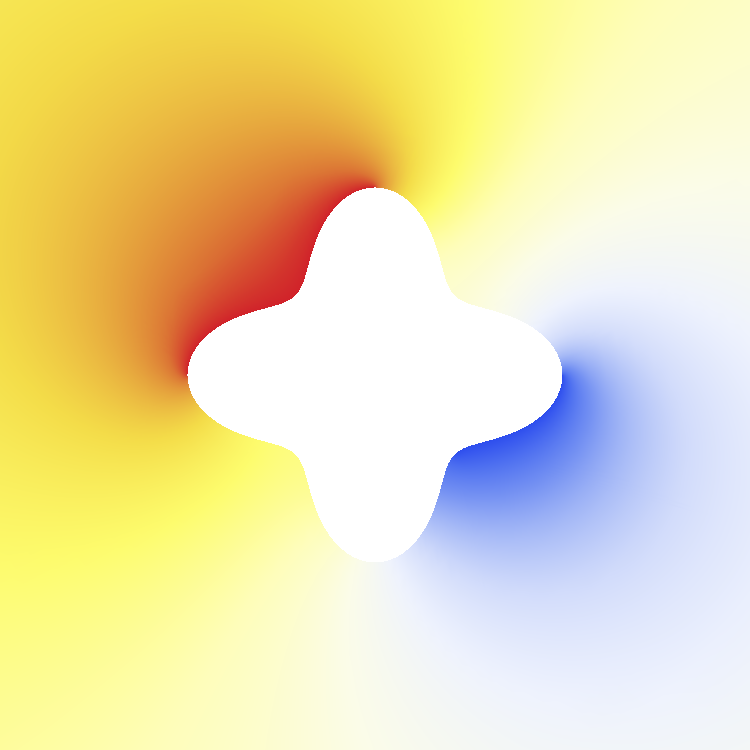}}
\caption{Quadrilateral C: The potential field of the inner (A) and outer (B) solutions.}\label{fig:inneroutercomparisonC}
\end{figure}

\begin{figure}
\centering
\subfloat[Potential field of the inner problem.]{\includegraphics[width=0.4\textwidth]{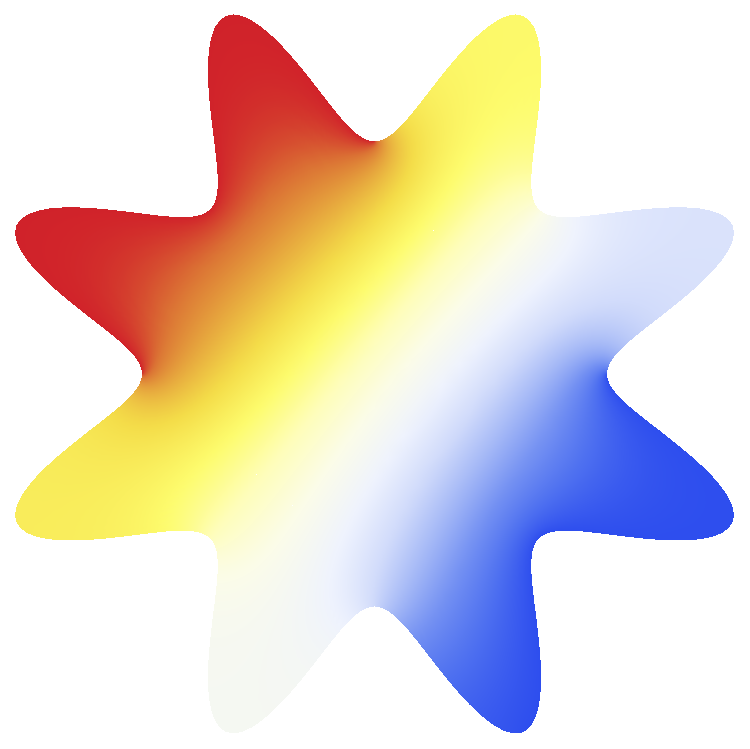}}
\subfloat[Potential field of the outer problem.]{\includegraphics[width=0.4\textwidth]{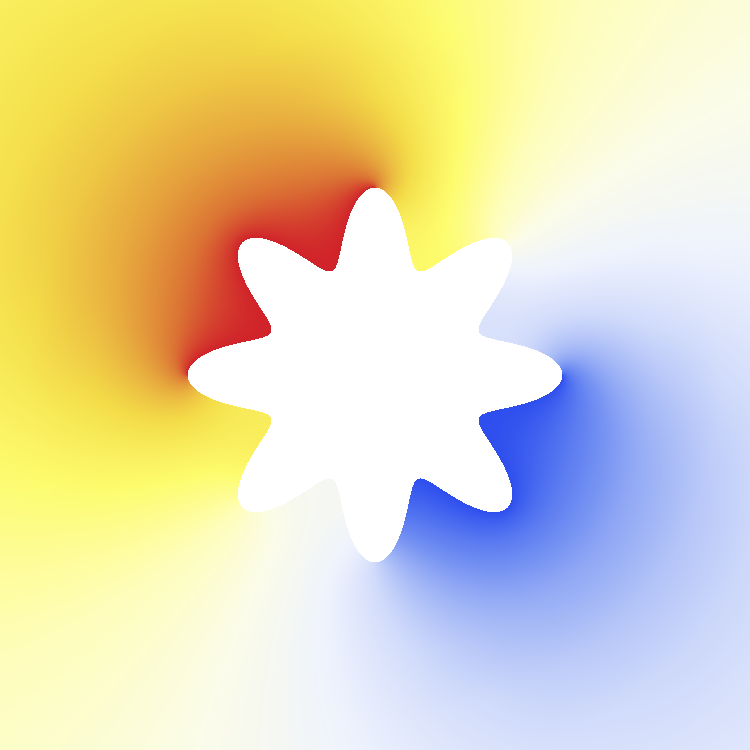}}
\caption{Quadrilateral D: The potential field of the inner (A) and outer (B) solutions.}\label{fig:inneroutercomparisonD}
\end{figure}

It turns out that besides the actual value of the exterior modulus one can also
determine the value of the far-field potential. Either one can determine the value
of the potential at the reflection point of the interior problem, i.e., at the origin,
or simply evaluate the solution of the exterior problem at the farthest point.
Remarkably, the truncated domain results agree well with the (theoretically) exact
results of the equivalent inner modulus problems (Table~\ref{tbl:infinitycomparison}).
In Figures~\ref{fig:inneroutercomparisonA}--\ref{fig:inneroutercomparisonD} we
show comparisons of the interior and exterior potential fields. For the two polygonal
quadrilaterals, the corresponding contour lines and the location of the origin
in the interior case are indicated.
In the general case, prediction of the far-field value based solely on geometric
arguments is an open problem.

We note, that for both Quadrilateral C and D, the interior and exterior capacities
are equal. This invariance is new and has not been reported in the literature before.
It is crucial that the four corners are chosen from extremal points, that is, local minima and maxima
of the radius.

\begin{table}
\begin{tabular}{|c|c|c|c|} \hline
Method & Exterior Modulus & Error (\ref{recipidty2})& Relative Error \\
\hline
SC Toolbox 			& 0.9923416323 	& -9 & -- \\
AFEM 				& 0.9923500126	& -4 &  -5 \\
$hp$-FEM (Interior) & 0.9923416332  & -9 &  -9 \\
$hp$-FEM (Exterior) & 0.9923416332  & -9 &  -9 \\
\hline
\end{tabular}
\caption{Quadrilateral A: $\{0, 1, ({28}/{25}, {69}/{50}), (-{19}/{25} , {21}/{25})\}$.
The values obtained with SC Toolbox are used as reference.
The errors are given as {$\lceil\log_{10}|\mathrm{error}|\rceil$}.}
\label{tbl:extmodcomparisonA}
\end{table}

\begin{table}
\begin{tabular}{|c|c|c|c|} \hline
Method & Exterior Modulus & Error (\ref{recipidty2})& Relative Error \\
\hline
SC Toolbox 			& 0.9592571721 	& -9 & -- \\
AFEM 				& 0.9593012739	& -4 &  -4 \\
$hp$-FEM (Interior) & 0.9592571731  & -8 &  -8 \\
$hp$-FEM (Exterior) & 0.9592572007  & -7 &  -7 \\
\hline
\end{tabular}
\caption{Quadrilateral B: $\{0,1, ({42}/{25}, 4), (-{3}/{25} , {21}/{25})\}$.
The values obtained with SC Toolbox are used as reference.
The errors are given as {$\lceil\log_{10}|\mathrm{error}|\rceil$}.}\label{tbl:extmodcomparisonB}
\end{table}

\begin{table}
\begin{tabular}{|c|c|c|c|} \hline
Method & Exterior Modulus & Error (\ref{recipidty2})& Relative Error \\
\hline
$hp$-FEM (Interior) & 0.8196441884805177  & -14 &  -- \\
$hp$-FEM (Exterior) & 0.8196441926483611  & -8 &  -8 \\
\hline
\end{tabular}
\caption{Quadrilateral C: $\gamma(t) = r(t) e^{it}, r(t) = 4/5 + (1/5) \cos(4 \pi t)$ and corners at $t=-1,-1/4,0,1/2$.
The values obtained with $hp$-FEM (Interior) are used as reference.
The errors are given as {$\lceil\log_{10}|\mathrm{error}|\rceil$}.}\label{tbl:extmodcomparisonC}
\end{table}

\begin{table}
\begin{tabular}{|c|c|c|c|} \hline
Method & Exterior Modulus & Error (\ref{recipidty2})& Relative Error \\
\hline
$hp$-FEM (Interior) & 0.9122187602015264  & -10 &  -- \\
$hp$-FEM (Exterior) & 0.9122187628550672  & -8 &  -8 \\
\hline
\end{tabular}
\caption{Quadrilateral D: $\gamma(t) = r(t)e^{it}, r(t) = 4/5 + (1/5) \cos(8 \pi t)$ and corners at $t=-1,-1/4,0,1/2$.
The values obtained with $hp$-FEM (Interior) are used as reference.
The errors are given as {$\lceil\log_{10}|\mathrm{error}|\rceil$}.}\label{tbl:extmodcomparisonD}
\end{table}

\begin{table}
\begin{tabular}{|c|c|c|c|} \hline
Quadrilateral & $hp$-FEM (Interior) & $hp$-FEM (Exterior) & Relative Error \\
\hline
A & 0.5281867366243582 	& 0.5281867468410989 & -7 \\
B & 0.6659476737428786	& 0.6659476800244547 & -8 \\
C & 0.5873283399651075  & 0.5873283469398137 & -7 \\
D & 0.5398927341965689  & 0.5398927414203410 & -7 \\ \hline
\end{tabular}
\caption{%
Comparison of the computed values of the potential at infinity.
The errors are given as {$\lceil\log_{10}|\mathrm{error}|\rceil$}.
}\label{tbl:infinitycomparison}
\end{table}

\section{Performance Considerations}

In this section we study the performance of our approach in terms of computational cost
in time and storage requirements, and convergence of the capacity, which is shown to be exponential.
Here we consider the Quadrilateral D defined above, and compare the interior and exterior problems.
This comparison is reasonable, since due to the new invariance, the interior and exterior problems
can be solved using \textit{exactly the same the geometry} and thus the singularities are of
the same kind.

\subsection{Convergence}
All experiments have been computed using $(\alpha, \nu)$-meshes, with $\alpha = 0.15$,
and $\nu = \min(16, p_{\max})$, where 16 is dictated by double precision. 
This choice allows us to compare two elemental $p$-distributions,
namely the constant $p = p_{\max}$, and the graded $p$-vector where the elemental $p$
increases per element layer away from the singularity, e.g., from $p=1$ up to $p = p_{\max}$.
The $p_{\max}$ has been chosen so that the relative error in both approaches is roughly the same
and in accordance with the results resported above, $p_{\max,I} = 18$ and  $p_{\max,E} = 22$, for the interior and exterior problems, respectively.

\begin{figure}
\centering
\subfloat[Loglog-scale: relative error in capacity vs $N$.]{\includegraphics[width=0.4\textwidth]{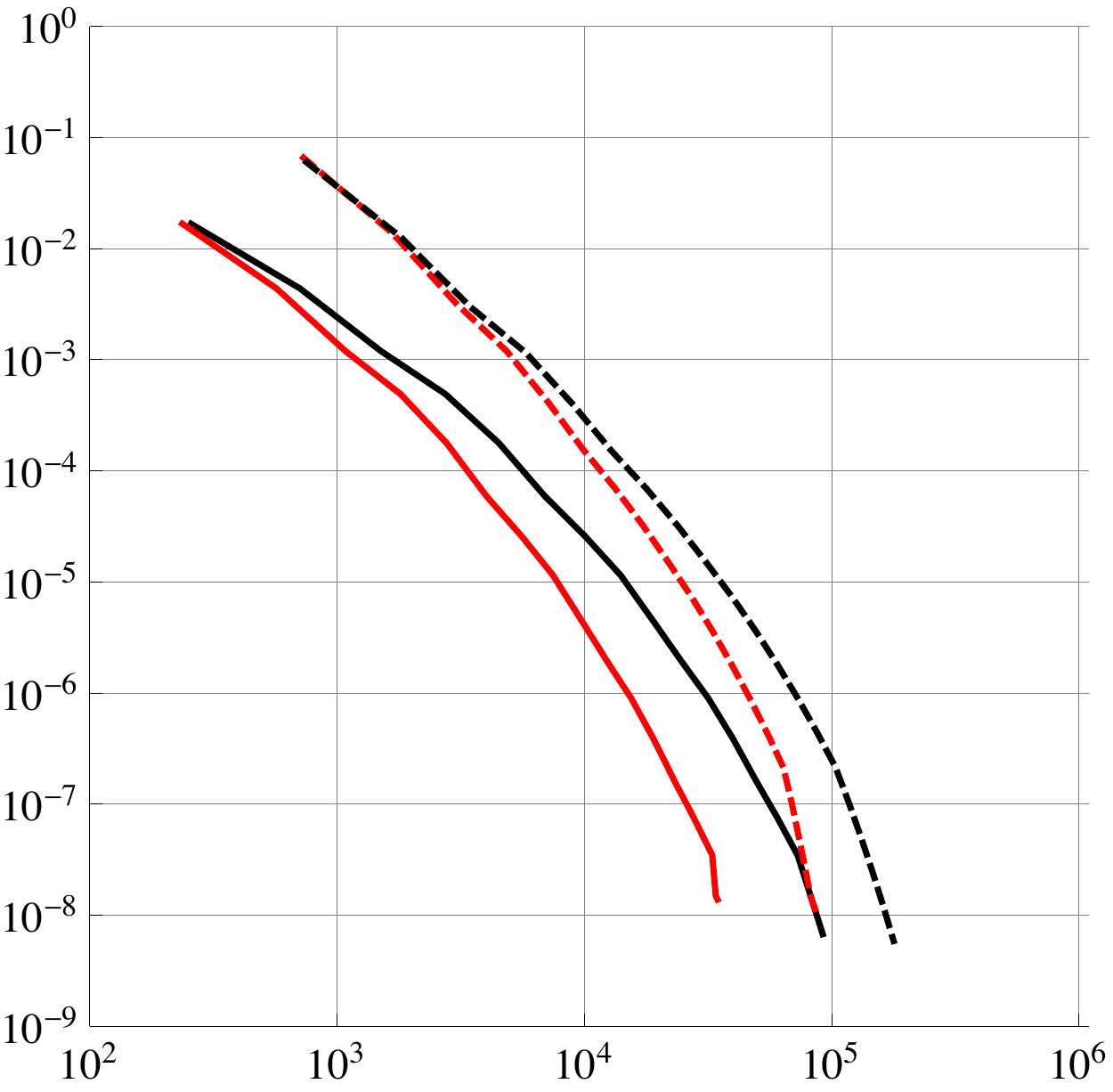}}\quad
\subfloat[Semilog-plot: relative error in capacity vs $N^{1/3}$.]
{\includegraphics[width=0.4\textwidth]{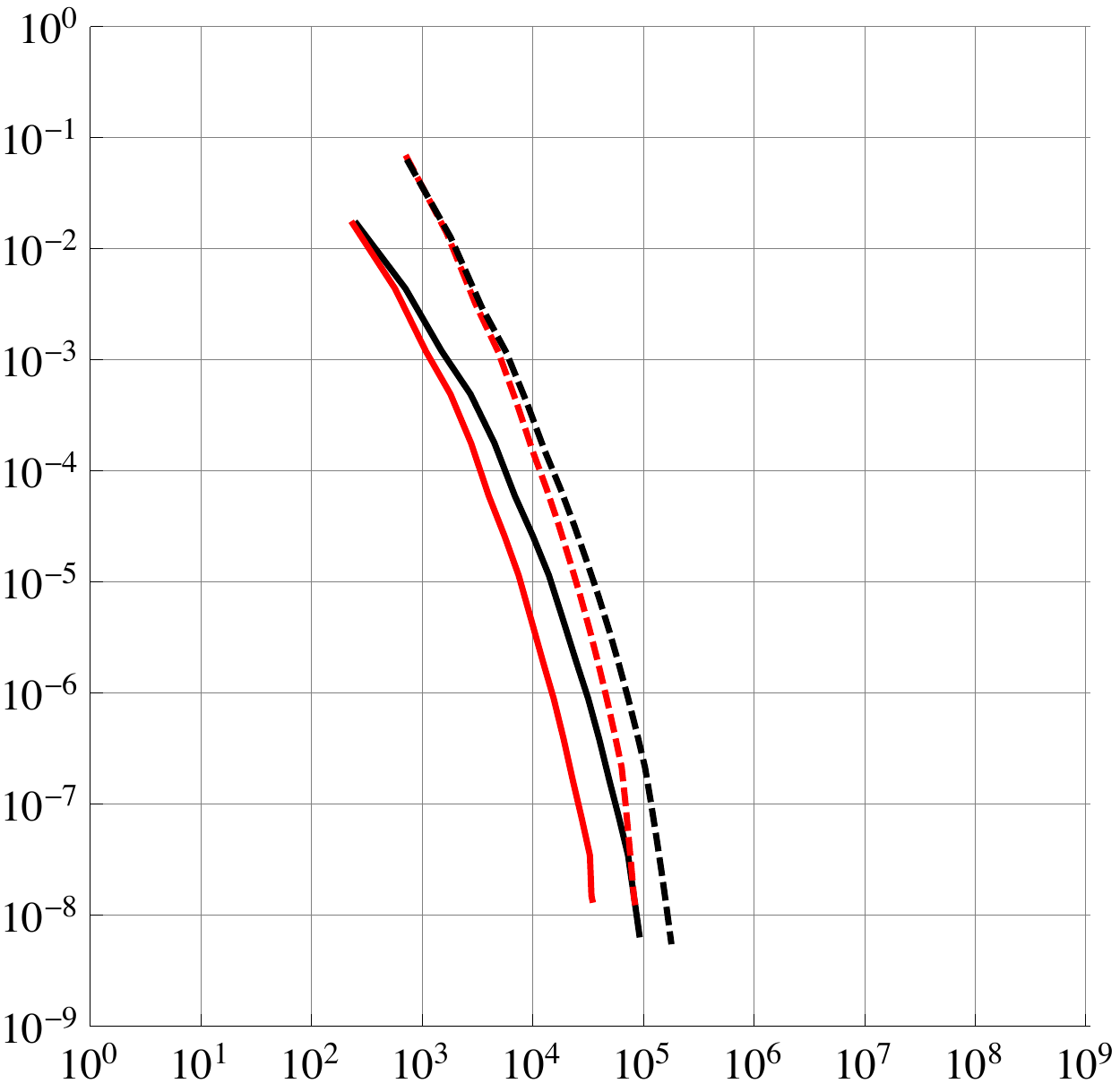}}
\caption{Quadrilateral D: Convergence of the relative error in capacity. Constant $p$: Black, Graded $p$: Red; Interior problem: Solid line; Exterior problem: Dashed line. The number of unknowns is $N$.}\label{fig:innerouterconvergence}
\end{figure}

The optimal rate of convergence of the relative error in capacity is 
\[\sim C_1\exp(-C_2\,N^{1/3}),\] 
where $N$ is the number of unknowns and $C_i$ are coefficients independent of $N$ \cite{s}.
In Figure~\ref{fig:innerouterconvergence} the convergence plots  
corresponding to both $p$-distributions
are shown using two different scalings: (A) in standard loglog-scale, and (B) in semilog-scale
with $N^{1/3}$ as the abscissa. 
The first plot shows that solutions to both problems converge exponentially,
but the latter one shows that the exterior approach is not as efficient as
the interior one.
Using linear fitting of logarithmic data, we find convergence rates of type $N^{1/\beta}$,
with $\beta_{I,c} = 3.72$, $\beta_{E,c} = 3.8$, $\beta_{I,g} = 3.41$, and $\beta_{E,g} = 3.55$,
where the indeces $c$ and $g$ refer to constant and graded polynomial distributions, respectively.

Two observations should be noted: a) faster convergence rate does not imply more accurate results;
b) the convergence behaviour becomes less stable as $p > \nu$ as the refinement strategy is
changed.

\begin{table}
\centering
\subfloat[Interior problem.]{
	\begin{tabular}{|c|c|c|rl|c|c|} \hline
		$p$ & $N$ & Meshing & Integration & (Assembly) & Solve & Total \\
		\hline
		4  & 1505  & 1  & 2   & (0)  & 0  & 3   \\
		8  & 10049 & 4  & 14  & (4)  & 4  & 22  \\
		12 & 31777 & 11 & 44  & (14) & 14 & 69  \\
		16 & 72833 & 21 & 136 & (52) & 42 & 199 \\
		\hline
	\end{tabular}
	}\\
\subfloat[Exterior problem.]{
	\begin{tabular}{|c|c|c|rl|c|c|} \hline
		$p$ & $N$ & Meshing & Integration & (Assembly) & Solve & Total \\
		\hline
		4  & 3456    & 4  & 2   & (0)  & 0  & 6   \\
		8  & 17792   & 13 & 19  & (6)  & 3  & 35  \\
		12 & 49152   & 26 & 63  & (23) & 10 & 99  \\
		16 & 103680  & 47 & 210 & (80) & 30 & 287 \\
		\hline
	\end{tabular}
	}
\caption{Quadrilateral D: Time spent in the solution process. All times are seconds as reported by Mathematica's Timing-function. Time spent in assembly of the linear system is included in that
of integration. (Apple Mac Pro 2009 Edition 2.26 GHz, Mathematica 8.0.4)}\label{tbl:inneroutercomparison}
\end{table}

\subsection{Time}
Averaged timing results over a set of 30 runs with constant $p$-distribution are shown in Table~\ref{tbl:inneroutercomparison}.
Note that the hierarchic nature of the problem has not been taken into account here and
runs for different values of $p$ have been independent. In our implementation the numerical integration
is the most expensive part. The numerical integration routines are based on a matrix-matrix multiplication
formalism which is highly efficient on terms of flops per memory access, and benefits from BLAS-level
parallelism on our test machine with eight cores; Apple Mac Pro 2009 Edition 2.26 GHz, Mathematica 8.0.4.
The time spent in assembling the matrix is included in the integration time. Mathematica does not
support pre-allocation of sparse matrix structures or autosumming initialization which leads to
a lot of reallocation of sparse matrices.

Interestingly, the time spent on direct solution of the systems is shorter for the exterior problem
for problems of comparable size.
In our opinion this is the result of the ordering heuristic used by Mathematica being more
efficient over ring domains.

\section{Conclusions}

In this study we have shown that three different algorithms, AFEM, SC Toolbox and $hp$-FEM, can all be effectively used for computation of the exterior modulus of a bounded polygonal quadrilateral. As far as we know, there are very few numerical or theoretical results on the exterior modulus in the literature.
The problem is first reduced to a Dirichlet-Neumann problem for the Laplace equation. In our earlier paper \cite{hrv} we introduced the reciprocal identity as an error estimate for the inner modulus computation
of a quadrilateral and here we demonstrate that the same method applies to error estimation for the exterior modulus as well.  We compare our numerical results to the analytic Duren-Pfaltzgraff formula for the exterior modulus of a rectangle and observe that our results agree with it. Moreover, in this case the analytic formula yields results that are within the limits provided by the reciprocal error estimate from our computational results. The reciprocal error estimate is also applied to study, for the case of polygonal quadrilaterals,
the accuracy of the Schwarz-Christoffel toolbox and the AFEM method, and mutual accuracy comparisons are given. Finally, for the case of quadrilaterals with curvilinear boundary, where these two methods do
not apply, we give  results obtained by the $hp$-FEM method, and their error estimates  based on the
relative error and the reciprocal error estimate. In this case we also analyze the dependence of the accuracy on the number of degrees of freedom and demonstrate nearly optimal convergence, compatible with
the theory of Babu\v ska and Guo \cite{bg1}.

A problem of independent interest is the value of the potential function at infinity. We study this problem for the exterior modulus of a polygonal quadrilateral and solve it by mapping the exterior domain onto a bounded domain by inversion and then computing the value of the potential function of the corresponding interior modulus problem at the image point of the point at infinity.

{\bf Acknowledgements.} The research of Matti Vuorinen was supported by the Academy of Finland, Project 2600066611. The authors are indebted to the referees for their helpful remarks.

\end{document}